\documentclass[review]{jsdarticle}
\usepackage{algorithm} \usepackage{algpseudocode} \usepackage{placeins} \usepackage{cleveref} \usepackage{tikz}
\usetikzlibrary{arrows.meta,positioning}

\curpubissue{For Review}{0}{2026}{1}{10.25518/2684-6500.xxx}
\curpubdates{2026-01-01}{2026-01-01}{2026-01-01}

\graphicspath{{./FIGURES/}}

\begin{document}

\title{Adjoint-Based Stochastic Analysis for Fatigue of Aeroelastic and Frictional Systems}

\author[1]{Justin H. Porter \email{Justin.Porter@nlr.gov}}

\affil[1]{National Laboratory of the Rockies, 15013 Denver West Parkway, Golden, CO, 80401, USA}

\begin{abstract}

As engineering systems across aerospace and energy generation applications are optimized to reduce weight and improve performance, vibration behavior, including nonlinear effects, is becoming increasingly important to analyze.
While methods such as harmonic balance are popular for analyzing periodic nonlinear vibration behavior, these analyses often only consider deterministic loading.
However, real applications experience additional noise in the excitation, and alternative approaches are frequently limited to smaller models.
This work expands on a recent approach to calculate a time-varying covariance using a linearized system around a periodic nonlinear response.
Specifically, numerical conditioning is addressed and improved for cases increasing from a few degrees of freedom to hundreds of degrees of freedom.
The algorithm is then extended to capture correlations across both space and time via a period-spanning covariance matrix.
Model reduction is proposed for this calculation to maintain computational efficiency.
The period-spanning covariance matrix is used to sample full stochastic cycles and estimate fatigue-relevant statistics up to 1393 times faster than stochastic time integration.
Importantly, the present methodology uses the same gradient information as stability analyses and continuation solutions, so it can easily be applied for many nonlinear analyses.
The methodology is demonstrated for two beam systems: one with unsteady aerodynamic forces and another with a Bouc-Wen frictional tip force.
For the aeroelastic model, statistics are accurately estimated via the present approach. 
For the frictional system, some errors are introduced compared to stochastic time integration and are attributed to the linear time-periodic nature of the system used in the present stochastic calculations. However, estimated stochastic statistics remain useful given the significant reduction in computation time.

\keywords{stochastic vibration; white noise; fatigue; aeroelasticity; friction; harmonic balance method; hysteresis}
\end{abstract}

\maketitle

\newcommand{\mat}[1]{\underline{\underline{#1}}}  \newcommand{\vect}[1]{\underline{#1}}  \renewcommand{\algorithmicrequire}{\textbf{Input:}}
\renewcommand{\algorithmicensure}{\textbf{Output:}}

\section{Introduction}

Vibration analysis is becoming increasingly important to ensure functionality and reliability as aerospace and energy generation systems become more highly optimized and flexible.
For instance, maximum expected vibration amplitudes are analyzed to ensure that unintended contact does not occur in jet engines and other applications.
Additionally, fatigue analysis is conducted to understand the expected useful life of structures.
Complicating analysis, engineering structures are becoming increasingly nonlinear with the inclusion of friction \cite{brakeMechanicsJointedStructures2017, berthold_fully_2023}, geometric \cite{krack_tribomechadynamics_2025, bagheriComputingPeriodicSE3_2025, touze_model_2021}, and aeroelastic \cite{dimitriadis_continuation_2008, branlard_using_2024, lee2020model, berthold_fully_2023} nonlinearities.
To calculate vibration amplitudes, computational models are used to simulate expected loading and vibration cases.
Such loadings can include regular periodic harmonic loads and more complex loads such as random/stochastic excitations or impacts.
The analysis of white noise excitations to a structure provides a useful starting point for analyzing more complicated stochastic loading cases
and for engineering design \cite{alotta_fractional_2014,branlard_using_2024, breunung_noise_2023, breunung_nonlinear_2025}.

Methods are well developed for analyzing vibration responses under deterministic loading.
First, models of structures are discretized (e.g., by the finite element method) and written as either a set of second-order differential equations or recast in state-space form as a set of first-order differential equations.
Time integration is frequently the most accessible simulation technique for nonlinear vibration, and variants of the generalized-$\alpha$ algorithm are commonly considered as benchmark algorithms \cite{bauchau_flexible_2011, bruls_lie_2012, wang_beamdyn_2017, shetty_fast_2020, berthold_fully_2023, witteveen_can_2023}.
However, time integration algorithms come with significant computational costs due to the sequential nature of the calculation.
Perturbation methods \cite{nayfehNonlinearOscillations1995,mathisReviewDampingModels2020, quinnIdentificationModeCoupling2023} and model reduction \cite{touze_model_2021, witteveenEfficientModeBased2009,kollerSemihyperReductionFiniteElement2020, cenedese_data_driven_2022, axas_fast_2023} offer opportunities to accelerate time integration, but remain an active area of research.
Alternatively, the harmonic balance method (HBM) is popular for directly calculating steady-state responses and is significantly faster than time integration \cite{krackHarmonicBalanceNonlinear2019, berthold_fully_2023}.
Likewise, shooting \cite{peeters_nonlinear_2009, peeters_theoretical_2010} and collocation \cite{karkar_comparative_2014, formica_coupling_2012, dankowicz_recipes_2013, ahsan_adjoint_based_2025} methods enable the direct calculation of deterministic steady-state responses.
Nonlinear modal analysis allows for rapid calculation of more varied responses \cite{krackNonlinearModalAnalysis2015, krackReducedOrderModeling2013, krack_computation_2014, raze_tracking_2024, porter_efficient_2025}.
Furthermore, uncertainty quantification has been considered for model parameters with polynomial chaos expansion \cite{yuan_propagation_2021, morsy_predicting_2025, balaji_vibration_2026} and stochastic HBM~\cite{didier_nonlinear_2013,sinou_stochastic_2015}. However, the efficient analysis of cases with stochastic excitation remains more difficult than deterministic analysis \cite{breunung_computationally_2022}.

For stochastic excitation, the It\^{o} stochastic differential equation is considered when using white noise excitation.
For traditional time integration methods, a sufficiently small time step can resolve the spectral content for some types of noise. For white noise, the energy content is spread across frequencies, so specialized time integration methods are necessary to resolve the responses.
While white noise may not be appropriate for many engineering applications, it can easily be adapted to other types of noise \cite{spanos_filter_1986, alotta_fractional_2014, breunung_noise_2023} and is a convenient basis for deriving analytical methods. Therefore, time integration schemes that can handle white noise are necessary to verify such analytical methods. 
The Euler-Maruyama method can be used for stochastic time integration \cite{kloeden_numerical_1999} but may require small time steps or long simulations to converge the statistics of the response.
Recently, a small noise integrator algorithm was developed that demonstrates computational speedups by applying better deterministic methods \cite{breunung_computationally_2022}.
However, the small noise integrator algorithm does not necessarily directly improve computation times when time steps for the deterministic portion of traditional methods (e.g., the generalized-$\alpha$ algorithm) drive the computation time.
Thus, alternatives to stochastic time integration are desirable.

\subsection{Fokker-Planck Equation}

Rather than calculating many stochastic trajectories, the distribution of the response can be directly computed with the Fokker-Planck-Kolmogorov equation or simply the Fokker-Planck equation (FPE) \cite{fuller_analysis_1969}.
While general solutions are difficult \cite{fuller_analysis_1969}, several numerical approaches have been proposed \cite{muscolino_stationary_1997,di_paola_approximate_2002,terejanu_uncertainty_2008}.
For related path integral formulations, computational costs can become prohibitive for systems with multiple degrees of freedom (DOFs), but modifications are sometimes possible \cite{kougioumtzoglou_efficient_2015, di_paola_path_2020}.
Alternatively, a number of simplifications to the FPE can be applied to ease solution computation under different assumptions, as outlined by \Cref{fig:fpe_simplifications}.

\FloatBarrier

\newcommand{\topleft}{FPE\\(Fully Nonlinear)\\\cite{fuller_analysis_1969,muscolino_stationary_1997,di_paola_approximate_2002,terejanu_uncertainty_2008}}
\newcommand{\topright}{Stationary FPE\\(Constant Equations, \\Nonlinear Stochastic)\\\cite{fuller_analysis_1969,muscolino_stationary_1997,di_paola_approximate_2002,kumar_multi_resolution_2006, kumar_partition_2006}}
\newcommand{\bottomleft}{Adjoint-based Covariance\\(Linear Time-Periodic or\\Quasiperiodic Stochastic) \cite{ahsan_adjoint_based_2025}} 
\newcommand{\bottomright}{Linear FPE\\(Fully Linear Time-Invariant)\\\cite{risken_fokker_planck_book_1996, kuehn_deterministic_2012, pavliotis_stochastic_2014,branlard_using_2024,kong_non_stationary_2022}}
\newcommand{\centerLeft}{Stochastic Averaging and \\ Statistical Linearization (Coupled \\ Nonlinear Deterministic and\\  Linear Stochastic) \cite{spanos_formulation_2019, spanos_nonlinear_2020, han_non_stationary_2022, kong_non_stationary_2022}} 
\newcommand{\labelTop}{Stationary/Time-Invariant}
\newcommand{\labelLeft}{Small\\(Linearized)\\Perturbations\\Near Limit\\Cycle/Torus}
\newcommand{\labelRight}{Small\\(Linearized)\\Perturbations}
\newcommand{\labelBottom}{Stationary/Time-Invariant}
\newcommand{\arrowToSA}{Perturbations \\on Average \\ Linearized System}
\newcommand{\arrowFromSA}{Fully Linear}

\newcommand{\height}[0]{7.5}
\newcommand{\width}[0]{10}

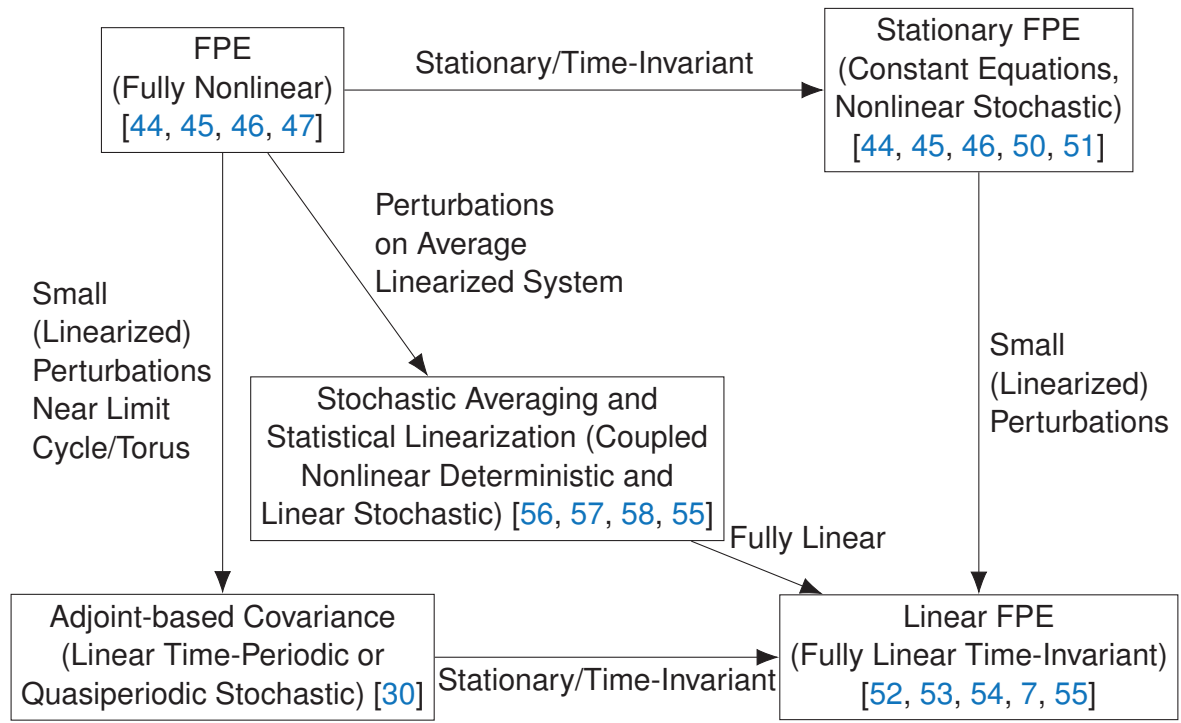
\begin{figure}[h!]
    \centering
        \begin{tikzpicture}[>=Latex, font=\sffamily, align=center]
        \node[draw, align=center] (A) at (0,\height) {\topleft};
    \node[draw, align=center] (B) at (\width,\height) {\topright};
    \node[draw, align=center] (C) at (0,0) {\bottomleft};
    \node[draw, align=center] (D) at (\width,0) {\bottomright};
    \node[draw, align=center] (E) at (0.35*\width,0.35*\height) {\centerLeft};
            \draw[-{Latex[length=3mm]}, align=center] (A) -- (B) node[midway,above] {\labelTop};
    \draw[-{Latex[length=3mm]}, align=left] (A) -- (C) node[midway,left] {\labelLeft};
    \draw[-{Latex[length=3mm]}, align=left] (B) -- (D) node[midway,right] {\labelRight};
    \draw[-{Latex[length=3mm]}, align=center] (C) -- (D) node[midway,below] {\labelBottom};
        \draw[-{Latex[length=3mm]}, align=left] (A) -- (E) node[pos=0.7, above right, xshift=-6pt] {\arrowToSA};
    \draw[-{Latex[length=3mm]}, align=center] (E) -- (D) node[pos=0.45, above right, xshift=-12pt] {\arrowFromSA};
        \end{tikzpicture}

    \caption{Diagram mapping fully nonlinear Fokker-Planck equation (FPE) to simpler alternatives.}
    \label{fig:fpe_simplifications}
\end{figure}

The largest simplification of the FPE is to a linear time-invariant system undergoing steady white noise excitation. This allows for the direct calculation of a covariance through the solution of a matrix problem (Lyapunov differential equation) \cite{risken_fokker_planck_book_1996, kuehn_deterministic_2012, pavliotis_stochastic_2014,branlard_using_2024, kong_non_stationary_2022}.
However, several intermediary simplifications are possible to provide more informative results.
Solutions can be found to the stationary FPE considering how perturbations respond nonlinearly \cite{fuller_analysis_1969,muscolino_stationary_1997,di_paola_approximate_2002,kumar_multi_resolution_2006, kumar_partition_2006}, but that approach is not directly applicable to cases that are nonstationary due to deterministic nonlinear responses.

Stochastic averaging and linearization provide an alternative for considering white noise excitation and either deterministic harmonic excitation \cite{anh_response_2014, spanos_formulation_2019, spanos_nonlinear_2020, kong_non_stationary_2022} or deterministic transient excitation \cite{han_non_stationary_2022}.
Stochastic averaging is used to transform the system response to slow variables (e.g., amplitude and phase) \cite{spanos_probability_1980, roberts_stochastic_1986, anh_response_2014, han_non_stationary_2022} and/or average the system description over a vibration cycle \cite{anh_response_2014, spanos_formulation_2019, spanos_nonlinear_2020}.
Additionally, a linearized description of the system can be used to calculate the response statistics due to the noise \cite{anh_response_2014,han_non_stationary_2022, spanos_formulation_2019, spanos_nonlinear_2020, kong_non_stationary_2022}.
The variance of the stochastic response, which is of interest here, influences the deterministic solution in some studies with two-way coupling effects \cite{spanos_formulation_2019, spanos_nonlinear_2020, han_non_stationary_2022, kong_non_stationary_2022}.
Specifically, a single harmonic HBM solution was coupled to the stochastic calculation in \cite{spanos_formulation_2019, spanos_nonlinear_2020}. In that work, the variance of the stochastic response influenced the deterministic solution, and the stochastic response was calculated based on a linearization around the deterministic solution.
Using stochastic averaging and linearization, studies have generally only considered small systems with 1 DOF \cite{spanos_probability_1980, roberts_stochastic_1986, han_non_stationary_2022, anh_response_2014} or 2 DOFs \cite{spanos_formulation_2019, kong_non_stationary_2022}.
Furthermore, the two-way coupling between HBM and the stochastic solution \cite{spanos_formulation_2019, spanos_nonlinear_2020} may introduce significant computational costs and implementation complexities for larger, engineering-scale systems.
Simultaneously, the averaging approach may eliminate oscillations in the variance that occur within a period of vibration \cite{han_non_stationary_2022}.

\FloatBarrier

\subsection{Adjoint Covariance Analysis}

Some of these limitations are addressed with the adjoint-based covariance approach \cite{ahsan_adjoint_based_2025}.
This approach allows for the covariance of the system to vary throughout a cycle by considering a time-varying linearization near a limit cycle or limit torus \cite{ahsan_adjoint_based_2025}. For a limit cycle, the stochastics are then calculated on a linear time-periodic (LTP) system.
This could enable the ability to capture oscillations in the variance that were averaged out in \cite{han_non_stationary_2022}.
Furthermore, the methodology proposed in \cite{ahsan_adjoint_based_2025} is a one-way coupling to calculate the stochastic response given the deterministic solution and uses derivative information that is generally already available in nonlinear vibration analyses for continuation and stability calculations \cite{kuehn_deterministic_2012, ahsan_adjoint_based_2025}. This likely reduces computational cost and eases implementation compared to the two-way coupled method of \cite{spanos_formulation_2019, spanos_nonlinear_2020} but may come with a reduction in accuracy. 
Given the benefits of the one-way coupled approach of \cite{ahsan_adjoint_based_2025}, some reduction in accuracy is acceptable for iterative engineering analysis phases.

While \cite{ahsan_adjoint_based_2025} develops a direct calculation of a periodic covariance matrix, further work is necessary to demonstrate the relevance of the approach to real engineering systems. 
Specifically, numerical stability can become a challenge for systems with many more DOFs.
Furthermore, a time-varying covariance matrix provides an incomplete understanding of the stochastic response where a response is expected to be partially correlated with responses at previous times.
For fatigue analysis, additional method development is required to understand the range and mean of vibration amplitudes over a cycle \cite{astm_fatigue_counting}.

The rest of this paper addresses these challenges with the adjoint-based stochastic calculation of \cite{ahsan_adjoint_based_2025}. First, the background theory is reviewed (\Cref{sec:background}), and then, updated and extended algorithms are presented in \Cref{sec:all_proposed_algs}.
Specifically, numerical conditioning is addressed in \Cref{sec:updated_alg}. To establish a more direct connection to fatigue analysis, Sections \ref{sec:period_spanning_cov} and \ref{sec:sampling_DPSC} develop a period-spanning covariance representation that allows for the direct sampling of full cycles of the stochastically excited response.
\Cref{sec:2dof_example} briefly reviews and extends the 2-DOF example from \cite{ahsan_adjoint_based_2025}.
Sections \ref{sec:aero_full} and \ref{sec:friction_system} present the approach applied to aeroelastic and frictional systems, respectively. These novel applications of the approach are highly relevant for real engineering design.
Lastly, conclusions are discussed in \Cref{sec:conclusions}.
Code for the present work and solution data are available at \cite{porterStochasticCodeData2026}.

\section{Background Theory} \label{sec:background}

This work analyzes the response of nonlinear systems to combined deterministic harmonic and white noise stochastic loading.
The system has second-order states $\vect{y}$ associated with the deformation of a beam and first-order states $\vect{z}$ that describe how internal variables for the nonlinear hysteretic (aerodynamic and frictional) forces evolve.
The second-order states have associated matrices for linear mass $\mat{M}$, damping $\mat{C}$, and stiffness $\mat{K}$. Additionally, nonlinear forces $\vect{f_{nl}}$ (detailed in Sections \ref{sec:aero_full} and \ref{sec:friction_system}) act on the second-order states and depend on both the first- and second-order states.
The full equations of motion are
\begin{equation} \label{eq:1stand2ndODE}
\begin{split}
    &\mat{M}\ \vect{\ddot{y}} + \mat{C}\ \vect{\dot{y}}  + \mat{K}\ \vect{y} + \vect{f_{nl}}(\vect{y}, \vect{\dot{y}}, \vect{z}, \vect{p}, t) = \vect{f_{ext}}(t)
    \\
    &\vect{\dot{z}} = \vect{g}(\vect{y}, \vect{\dot{y}}, \vect{z}, \vect{p}, t),
\end{split}
\end{equation}
where a dot indicates a time derivative. Additionally, a time-varying external force $\vect{f_{ext}}$ is applied to the second-order states, and the derivative of first-order states is a general nonlinear function $\vect{g}$.
The nonlinear forces $\vect{f_{nl}}$ and the function $\vect{g}$ may both depend on some system parameters in the vector $\vect{p}$ and the time $t$.

For the remainder of this paper, the methodology of \cite{ahsan_adjoint_based_2025} in the case of a limit cycle is referred to as the direct periodic covariance (DPC).
Since DPC \cite{ahsan_adjoint_based_2025} considers only first-order systems, the states are combined into a vector of $N$ first-order states $\vect{x}$ as 
\begin{equation}
    \vect{x} = 
    \begin{bmatrix}
        \vect{y} \\
        \vect{\dot{y}} \\
        \vect{z}
    \end{bmatrix}.
\end{equation}
The response of the nonlinear system is described by rearranging \eqref{eq:1stand2ndODE} into a state-space form giving the vector-function $\vect{f}$:
\begin{equation}
    \vect{\dot{x}} = \vect{f}(t, \vect{x}, \vect{p}).
\end{equation}
Additionally, this work considers single harmonic deterministic excitation.
Time integration for deterministic and stochastic excitation and HBM are described in Appendix~\ref{sec:appendix_methods}.

As the basis of DPC, the derivative matrix $\mat{X}(t)$ is defined as
\begin{equation}
    \mat{X}(t) = \dfrac{\partial \vect{x}(t)}{\partial \vect{x_0}}.
\end{equation}
This is commonly calculated for Floquet analysis (e.g., see \cite{peeters_theoretical_2010, peeters_nonlinear_2009}), which uses the monodromy matrix of $\mat{X}(T)$ for period of vibration $T$.
Using the first-order state representation of the system $\vect{x}$, the derivative propagation over a time step\footnote{Specifically, this is derived for the trapezoid time integration scheme.} of size $\Delta t$ is 
\begin{equation}
    \dfrac{\partial \vect{x}^{n+1}}{\partial \vect{x}^n} 
    = 
    - \left[ \dfrac{\Delta t}{2}\mat{A}^{n+1} - \mat{I_N}  \right]^{-1} \left[ \dfrac{\Delta t}{2} \mat{A}^n+ \mat{I_N} \right].
\end{equation}
For time index denoted with superscript $j$,
\begin{equation}
    \mat{A}^j = \dfrac{\partial \vect{f}(t, \vect{x}, \vect{p})}{\partial \vect{x}} \Bigg|_{\vect{x} = \vect{x}^j, t = t^j},
\end{equation}
and $\mat{I_N}$ is the $N\times N$ identity matrix.
A discrete-time representation of $\mat{X}(t)$ can be calculated via the chain rule using the derivatives over individual time steps and with initial conditions of $\mat{X}(0) = \mat{I_N}$.

\subsection{Original Direct Periodic Covariance}

The original DPC method of \Cref{alg:original_direct} provides a direct computation for the time-varying covariance of a steady-state periodic system under white noise excitation \cite{ahsan_adjoint_based_2025}.
Only the periodic case is considered here, though \cite{ahsan_adjoint_based_2025} also addresses quasi-periodic deterministic excitation plus noise.
\Cref{alg:original_direct} is based on a cycle of vibration of length $T$ seconds divided into $N_t$ time steps and returns the covariance matrix $\mat{C}^j$ at each time step.
The algorithm takes inputs of $\mat{X}^j$, which is the discrete form of $\mat{X}(t^j)$, and $\mat{F}^j$, which describes how the stochastic forcing is applied.
In this work, parameters used in the nonlinear force calculation are assumed to directly undergo white noise excitation. 
Throughout this work, the number of time steps per period, $N_t$, is assumed to be even for formulating algorithms and is always used as a power of 2, consistent with adopting the alternating frequency-time (AFT) algorithm for HBM (see \Cref{sec:appendix_methods_hbm}).

For white noise variations in the parameter $p_i$, the stochastic forcing vector is calculated as
\begin{equation}
    \mat{F}(t, \vect{x}, \vect{p}) = \dfrac{\partial \vect{f}(t, \vect{x}, \vect{p})}{\partial p_i}.
\end{equation}
This work uses the derivative only with respect to a single parameter per case, so $\mat{F}$ is of size $ N \times 1$ and describes how the stochastic deviation effects get distributed to the system DOFs.
Note that \cite{ahsan_adjoint_based_2025} utilizes a scaling of $\sigma_0\sqrt{T}$ for the random forcing vector $\mat{F}^j$. Here, this scaling is incorporated into the full standard deviation of the white noise as $\sigma = \sigma_0 \sqrt{T}$ to simplify the notation.
The result $\mat{C}^j$ of \Cref{alg:original_direct} can be scaled by $\sigma^2$ for a chosen value of $\sigma$ to get the desired covariance.

\FloatBarrier

\newcommand{\indicesNotation}[0]{Indices are denoted by superscripts, start at zero, and include endpoints.}

\begin{algorithm}[h!]
\caption{Original calculation of the DPC as presented by \cite{ahsan_adjoint_based_2025} and using a numerical integral for $\mathcal{I}^j$. Resulting $\mat{C}^j$ is a time-periodic covariance and can be directly scaled by the noise variance $\sigma^2$. \indicesNotation}\label{alg:original_direct}
\begin{algorithmic}[1]
\Require $N_t$, $T$, $\mat{X}^j$, $ \mat{F}^j$ $\forall j \in [0, N_t]$
\State $dt \gets T / N_t$
\For{$j=0$ to $N_t$}
    \State $\mat{G}^j \gets \left[\mat{X}^j \right]^{-1} \mat{F}^j $  \label{alg:alg1_stepXinv}
    \State $\mat{d\mathcal{I}}^j \gets \mat{G}^j \left[  \mat{G}^j \right]^T$
\EndFor
\State $\mat{\mathcal{I}}^0 \gets \mat{0_N}$ \Comment{$\mat{0_N}$ is an $N \times N$ zero matrix.}
\For{$j = 2$ to $N_t$ with step 2}
    \State $\mat{\mathcal{I}}^{j-1} \gets \dfrac{dt}{12} \left[5 \mat{d\mathcal{I}}^{j-2} + 8 \mat{d\mathcal{I}}^{j-1 } - \mat{d\mathcal{I}}^{j} \right] $
    \State $\mat{\mathcal{I}}^j \gets \dfrac{dt}{3} \left[\mat{d\mathcal{I}}^{j-2} + 4 \mat{d\mathcal{I}}^{j-1 } + \mat{d\mathcal{I}}^{j} \right] $
\EndFor
\State $\mat{\mathcal{J}}^{N_t} \gets \mat{X}^{N_t} \left[ \mat{\mathcal{I}}^{N_t} \right]  \left[ \mat{X}^{N_t} \right]^T $
\State $\mat{C_0} \gets \text{discreteLyapunovSolution}(\mat{X}^{N_t}, \mat{\mathcal{J}}^{N_t})$ \Comment{Solution of $\mat{C_0} = \mat{X}^{N_t} \mat{C_0} \left[ \mat{X}^{N_t} \right]^T + \mat{\mathcal{J}}^{N_t}$}
\label{alg:alg1_stepDiscreteLyap}
\For{$j=0$ to $N_t$}
    \State $\mat{C}^j \gets \mat{X}^j \mat{C_0} \left[ \mat{X}^j \right]^T + \mat{X}^j \mat{\mathcal{I}}^j \left[ \mat{X}^j \right]^T$
\EndFor
\Ensure $\mat{C}^j \ \forall j \in [0, N_t]$
\end{algorithmic}
\end{algorithm}

The solution of the discrete Lyapunov equation in Step \ref{alg:alg1_stepDiscreteLyap} of \Cref{alg:original_direct} can be obtained directly with the calculations
\begin{equation} \label{eq:directLyapunov}
\begin{split}
    \mat{M} &= \mat{I_{N2}} - \mat{X}^{N_t} \otimes  \mat{X}^{N_t}  \\
        \mat{C_0} &= \text{matrix}\left( \mat{M}^{-1} \text{vector}\left( \mat{\mathcal{J}}^{N_t} \right) \right),
\end{split}
\end{equation}
where $\mat{I_{N2}}$ is the $N^2 \times N^2$ identity matrix.
However, for large systems this approach becomes computationally prohibitive, so faster algorithms are available \cite{barraud_numerical_1977}.

\FloatBarrier

\section{Proposed Algorithms} \label{sec:all_proposed_algs}

This section provides the details for the algorithms used throughout this paper.
First, \Cref{sec:updated_alg} discusses modifications to the DPC to improve numerical conditioning. The approach is then extended to include the covariance across different times in \Cref{sec:period_spanning_cov}. Lastly, the sampling approach to efficiently estimate fatigue-relevant statistics is presented in \Cref{sec:sampling_DPSC}.

\subsection{Modified Direct Periodic Covariance} \label{sec:updated_alg}

This work focuses on the periodically excited limit cycle case presented in \cite{ahsan_adjoint_based_2025}. 
The quasi-periodic limit tori case is not considered, and self-excited systems are briefly discussed in Appendix~\ref{sec:app_conditioned_self_excited}.
The present algorithm for the DPC is outlined in \Cref{alg:conditioned_direct}. A key change is the calculation of $\mat{\mathcal{J}}$, which is discussed below.
Since this work considers systems with more DOFs, the direct method for the discrete Lyapunov equation (see \eqref{eq:directLyapunov}) is not practical, so more efficient methods should be used \cite{barraud_numerical_1977}. Specifically, this work uses the implementation in Python's SciPy.\footnote{\url{https://docs.scipy.org/doc/scipy/reference/generated/scipy.linalg.solve_discrete_lyapunov.html}}

\FloatBarrier

\begin{algorithm}[h!]
\caption{Modified algorithm for DPC that improves numerical conditioning. \indicesNotation}\label{alg:conditioned_direct}
\begin{algorithmic}[1]
\Require $N$, $N_t$, $T$, $\mat{A}^j$, $ \mat{F}^j$ $\forall j \in [0, N_t]$
\State $dt \gets T / N_t$
\State $\mat{X}^j,  \mat{\mathcal{J}}^j \gets \text{integrateXJ}\left(dt, N_t, \mat{A}^j, \mat{F}^j \ \forall j \in [0, N_t] \right) $ \Comment{Evaluated with \Cref{alg:integrateJ}.}
\State $\mat{C_0} \gets \text{discreteLyapunovSolution}(\mat{X}^{N_t}, \mat{\mathcal{J}}^{N_t})$ \Comment{Solution of $\mat{C_0} = \mat{X}^{N_t} \mat{C_0} \left[ \mat{X}^{N_t} \right]^T + \mat{\mathcal{J}}^{N_t}$}
\For{$j=0$ to $N_t$}
    \State $\mat{C}^j \gets \mat{X}^j \mat{C_0} \left[ \mat{X}^j \right]^T + \mat{\mathcal{J}}^j $
\EndFor
\Ensure $\mat{X}^j,  \mat{C}^j \ \forall j \in [0, N_t]$
\end{algorithmic}
\end{algorithm}

In \Cref{alg:original_direct}, the calculation of $\mat{G}^j$ is likely to be ill-conditioned for numerical engineering models. For $j = N_t$, the calculation of $\mat{G}^{N_t}$ requires inverting the monodromy matrix $\mat{X}^{N_t} = \mat{X}(T)$.
For an engineering model with $N$ DOFs, the lower modes are generally of the greatest interest, and the highest modes become numerical artifacts resulting from the finite element discretization. To obtain useful solutions, these numerical higher modes will generally be positively damped and may have significant damping in the case of stiffness-proportional damping models. Since the period of these modes is generally much smaller than the vibration period of interest, $T$, the high-frequency modes will undergo many vibration cycles during one period of interest and thus significantly decay. Thus, the matrix $\mat{X}^{N_t}$ is likely to have near zero eigenvalues, which in double precision may be numerically zero. 
One could use an approximate matrix solve in Step \ref{alg:alg1_stepXinv} in \Cref{alg:original_direct}, but later examples demonstrate that this is not consistently sufficient to ensure a good final result (see Sections \ref{sec:aero_conditioning} and \ref{sec:fric_conditioning}).

To avoid this issue, the calculation of $\mat{\mathcal{J}}^j$ is rearranged in this work, as presented in \Cref{alg:integrateJ}. 
The key observation is that $\mat{G}^i$ is never used in isolation; rather, it is always pre-multiplied by $\mat{X}^j$ for $j \geq i$. 
In the continuous form \cite{ahsan_adjoint_based_2025}, $\mat{\mathcal{J}}$ is
\begin{equation}
\begin{split}
    \mat{\mathcal{J}}(t^j) 
    &= \mat{X}(t^j) \int_0^{t^j} \mat{G}(s) \mat{G}^T(s) ds \mat{X}^T({t^j})
    \\ &=  \int_0^{t^j} \left[ \mat{X}({t^j}) \mat{X}^{-1}(s) \right] \mat{F}(s) \mat{F}^T(s) \mat{X}^{-T}(s) \mat{X}^T({t^j}) ds. 
\end{split}
\end{equation}
The bracketed quantity for $s = t^i$ is
\begin{equation}
     \mat{X}({t^j}) \mat{X}^{-1}(t^i) = 
     \dfrac{\partial \vect{x}(t^j)}{\partial \vect{x}(t^i)}.
\end{equation}
This quantity is a segment of the integral of the monodromy matrix and thus can be directly calculated using the same approach as the monodromy matrix without requiring the inversion of the likely ill-conditioned $\mat{X}(t^i)$.
\Cref{alg:integrateJ} presents the updated algorithm for the discrete evaluation of $\mat{\mathcal{J}}$. Since \Cref{alg:integrateJ} requires much of the intermediate information to evaluate $\mat{X}^j$, it is written to include the computation of $\mat{X}^j$ rather than requiring $\mat{X}^j $ as input to \Cref{alg:conditioned_direct}.
\Cref{alg:integrateJ} uses Simpson's rule for the even steps to get better accuracy and the trapezoid rule to fill in the data at the odd time steps. Note that the evaluations of $\mat{\mathcal{J}}^j$ for odd $j$ do not directly influence the subsequent integral and only affect the computation of $\mat{C}^j$ for the same $j$.

\begin{algorithm}[h!]
\caption{Algorithm for integration of $\mat{\mathcal{J}}$ and $\mat{X}$. \indicesNotation}\label{alg:integrateJ}
\begin{algorithmic}[1]
\Require $N$, $N_t$, $T$, $\mat{A}^j$, $\mat{F}^j$ $\forall j \in [0, N_t]$
\State $dt \gets T / N_t$
\State $\mat{\mathcal{J}}^0 \gets \mat{0_N}$ \Comment{$\mat{0_N}$ is an $N \times N$ zero matrix.}
\State $\mat{X}^0 \gets \mat{I_N}$ \Comment{$\mat{I_{N}}$ is the $N \times N$ identity matrix.}
\For{$j = 0$ to $(N_t-2)$ with step 2}

    \State $\mat{dX}^{j+1} \gets \left[ \mat{I_{N}} - \dfrac{dt}{2} \mat{A}^{j+1} \right]^{-1} \left[ \mat{I_{N}} + \dfrac{dt}{2} \mat{A}^{j} \right] $
    \State $\mat{dX}^{j+2} \gets \left[ \mat{I_{N}} - \dfrac{dt}{2} \mat{A}^{j+2} \right]^{-1} \left[ \mat{I_{N}} + \dfrac{dt}{2} \mat{A}^{j+1} \right] $

    \State $\mat{d\mathcal{J}}^{j} \gets \mat{F}^{j} \left[ \mat{F}^{j} \right]^T $
    \State $\mat{d\mathcal{J}}^{j+1} \gets \mat{F}^{j+1} \left[ \mat{F}^{j+1} \right]^T $
    \State $\mat{d\mathcal{J}}^{j+2} \gets \mat{F}^{j+2} \left[ \mat{F}^{j+2} \right]^T $

    \State $\mat{\tilde{\mathcal{J}}}^{j+1} \gets \mat{dX}^{j+1} \mat{\mathcal{J}}^{j} \left[ \mat{dX}^{j+1} \right]^T $ 
    \Comment{Advance previous portion of the integral.}
    
    \State $\mat{\tilde{\mathcal{J}}}^{j+2} \gets \mat{dX}^{j+2} \mat{dX}^{j+1} \mat{\mathcal{J}}^{j} \left[ \mat{dX}^{j+1} \right]^T  \left[ \mat{dX}^{j+2} \right]^T $

    \State $\mat{{\mathcal{J}}}^{j+1} \gets \mat{\tilde{\mathcal{J}}}^{j+1} + \dfrac{dt}{2} \left[ 
    \mat{dX}^{j+1} \mat{d\mathcal{J}}^{j} \left[ \mat{dX}^{j+1} \right]^T
    +
    \mat{d\mathcal{J}}^{j+1} 
    \right] $
    \Comment{Add integral from $j$ to $j+1$}
    
    \State $\mat{{\mathcal{J}}}^{j+2} \gets \mat{\tilde{\mathcal{J}}}^{j+2} 
    + \dfrac{dt}{3} \left[ 
    \mat{dX}^{j+2} \mat{dX}^{j+1} \mat{d\mathcal{J}}^{j} \left[ \mat{dX}^{j+2} \mat{dX}^{j+1}  \right]^T
    + 4 \mat{dX}^{j+2} \mat{d\mathcal{J}}^{j+1} \left[ \mat{dX}^{j+2} \right]^T
    + \mat{d\mathcal{J}}^{j+2} 
    \right] $

    \State $\mat{X}^{j+1} \gets \mat{dX}^{j+1} \mat{X}^j $ \Comment{Progress integration of monodromy matrix.}
    \State $\mat{X}^{j+2} \gets \mat{dX}^{j+2} \mat{X}^{j+1} $

\EndFor
\Ensure $\mat{X}^j,  \mat{\mathcal{J}}^j \ \forall j \in [0, N_t]$
\end{algorithmic}
\end{algorithm}

This work does not consider self-excited systems. However, a similar modification to \cite{ahsan_adjoint_based_2025} is possible to the computation of $\mat{\mathcal{J}}$ in that case to avoid the ill-conditioned matrix solves.
Since the self-excited case is not used here, the additional modification is only outlined in Appendix~\ref{sec:app_conditioned_self_excited}.

\FloatBarrier

\subsection{Direct Period-Spanning Covariance} \label{sec:period_spanning_cov}

This section develops a direct period-spanning covariance (DPSC) that can be directly calculated without iteration or random sampling and describes covariance relations between both spatial DOFs and discrete times.
For fatigue analysis, it is desired to estimate the distribution of the range between the maximum and minimum amplitude during a vibration cycle without prescribing the time at which those occur.
Specifically, one cannot guarantee the time at which the maximum vibration amplitude occurs, so the time-varying covariance calculated in \Cref{sec:updated_alg} is insufficient for estimating fatigue-relevant statistics.
Furthermore, the variance at one point in the cycle cannot be treated as independent of the variance at a later point. 
The covariance between two random vectors $\vect{x}^i$ and $\vect{x}^j$ is
\begin{equation}
    \mat{C}^{ij} = E\left[ 
    (\vect{x}^i - \vect{\bar{x}}^i)(\vect{x}^j - \vect{\bar{x}}^j)^T \right],
\end{equation}
where $E \left[\cdot\right]$ denotes the statistical expectation and $\bar{\cdot}$ denotes the mean.
When $ i = j$, this reduces to $\mat{C}^j$, as calculated in \Cref{sec:updated_alg}.
For $j > i$, the perturbation away from the mean solution at time $j$ is directly related to the perturbation at time $i$ as (e.g., see \cite{risken_fokker_planck_book_1996})
\begin{equation}
    \vect{x}^j - \vect{\bar{x}}^j = 
        \dfrac{\partial \vect{x}(t^j)}{\partial \vect{x}(t^i)}
    \left[ \vect{x}^i - \vect{\bar{x}}^i \right]
    +
    \vect{x_{rand}}(t^i, t^j),
\end{equation}
where the first term on the right describes how the perturbation (generally) decays with time and the second term describes the effect of the white noise excitation on the system between times $t^i$ and $t^j$. 
It is noted that since the white noise excitation is uncorrelated between times, $\vect{x_{rand}}(t^i, t^j)$ is uncorrelated from $\vect{x}^i - \vect{\bar{x}}^i$. Thus,
\begin{equation}
\begin{split}
    \mat{C}^{ij} &= E\left[ 
    (\vect{x}^i - \vect{\bar{x}}^i)\left( \dfrac{\partial \vect{x}(t^j)}{\partial \vect{x}(t^i)}
    \left[ \vect{x}^i - \vect{\bar{x}}^i \right]
    +
    \vect{x_{rand}}(t^i, t^j) \right)^T \right]
    \\
    &= E\left[ 
    (\vect{x}^i - \vect{\bar{x}}^i)
    \left[ \vect{x}^i - \vect{\bar{x}}^i \right]^T 
    \left[ \dfrac{\partial \vect{x}(t^j)}{\partial \vect{x}(t^i)} \right]^T \right]
    \\
    &= \mat{C}^i \left[ \dfrac{\partial \vect{x}(t^j)}{\partial \vect{x}(t^i)} \right]^T,
\end{split}
\end{equation}
where the partial derivative term can be moved outside of the statistical expectation under the assumption of an LTP system, as is also employed for the derivation of the DPC.
Using the same discrete representation as \Cref{sec:updated_alg} for the partial derivative, \Cref{alg:period_spanning} presents a method for calculating a single DPSC matrix that captures both spatial and temporal correlations between DOFs. Here, the DPSC matrix is calculated for temporal correlations over a single period, but it could easily be extended to cover more time at the cost of a larger matrix.

\FloatBarrier

\begin{algorithm}[h!]
\caption{Algorithm for calculating DPSC, which describes covariance between times across a period. \indicesNotation}\label{alg:period_spanning}
\begin{algorithmic}[1]
\Require $N$, $N_t$, $T$, $\mat{C}^j$, $ \mat{A}^j$ $\forall j \in [0, N_t]$
\State $\mat{C_P} \gets$ zeros($(N_t + 1)N$, $(N_t + 1)N$)
\State $dt \gets T / N_t$
\For{$j=0$ to $N_t$}
\State $J \gets j N $ to $ (j+1) N - 1 $ \Comment{$J$ refers to the $j$th matrix block of size $N$.}
\State $ \mat{C_P}^{J,J} \gets \mat{C}^j$ \Comment{Implementation enforces symmetry here.}
\If{$j > 0$}
\State $(J-1) \gets (j-1) N $ to $ j N - 1 $ \Comment{$(J-1)$ refers to the $(j-1)$th matrix block.}
\State $ \mat{dX}^j \gets \left[ \mat{I_{N}} - \dfrac{dt}{2} \mat{A}^j \right]^{-1} \left[ \mat{I_{N}} + \dfrac{dt}{2} \mat{A}^{j-1} \right]$
\Comment{$\mat{I_{N}}$ is the $N \times N$ identity matrix.}
    \For{$k=0$ to $j - 1$}
            \State $K \gets k N $ to $ (k+1) N - 1 $ \Comment{$K$ refers to the $k$th matrix block.}
                \State $ \mat{C_P}^{K,J} \gets \mat{C_P}^{K,J-1} \left[\mat{dX}^j\right]^T$
        \State $ \mat{C_P}^{J,K} \gets \left[ \mat{C_P}^{K,J} \right]^T$
        \EndFor 
\EndIf
\EndFor
\Ensure $\mat{C_P}$
\end{algorithmic}
\end{algorithm}

\Cref{alg:period_spanning} generates a matrix of size $N(N_t+1) \times N (N_t + 1)$. Given that a value of $N_t$ of 128 is common for AFT calculations (and thus provides a reasonable starting point here), this matrix could become much larger than $N\times N$.
Rather, it is desirable to use a reduced basis. In \Cref{alg:cov_across_time_rom}, this work proposes a reduced basis of size $m$ at each time corresponding to the covariance shapes, with the largest covariance at each time step (i.e., the eigenvectors\footnote{Since $\mat{C}^j$ is symmetric and at least positive semi-definite, this basis is equivalent to one calculated via a singular value decomposition.} corresponding to the $m$ largest eigenvalues of $\mat{C}^j$).
Furthermore, eigenvectors $\mat{V}^j$ of $\mat{C}^j$ are returned from the algorithm to allow for the recovery of sampled statistics to the original coordinates.

\FloatBarrier

\begin{algorithm}[h!]
\caption{Algorithm for reduced-order model of covariance between times across a period. \indicesNotation}\label{alg:cov_across_time_rom}
\begin{algorithmic}[1]
\Require $N_t$, $m$, $T$, $\mat{C}^j$, $ \mat{A}^j$ $\forall j \in [0, N_t]$
\State $\mat{C_P} \gets$ zeros($(N_t + 1)m$, $(N_t + 1)m$) \Comment{Approximate with $m$ vectors at each time.}
\State $dt \gets T / N_t$
\For{$j=0$ to $N_t$}
\State $\vect{\lambda}^j, \mat{V}^j \gets \text{eigs}\left(\mat{C}^j, m, \text{`largest'}\right)$ \Comment{Calculate $m$ eigenpairs with largest eigenvalues.}
\State $\mat{\Lambda}^j \gets \text{diag}(\vect{\lambda}^j)$ \Comment{Convert eigenvalues to a diagonal matrix.}
\State $J \gets j m $ to $ (j+1) m - 1 $ \Comment{$J$ refers to the $j$th matrix block of size $m$.}
\State $ \mat{C_P}^{J,J} \gets \mat{\Lambda}^j $
\If{$j > 0$}
        \State $ K \gets 0 $ to $ j m - 1 $ \Comment{$K$ refers to all of the first $(j-1)$ matrix blocks.}
        \State $ \mat{dX}^j \gets \left[ \mat{I_{N}} - \dfrac{dt}{2} \mat{A}^j \right]^{-1} \left[ \mat{I_{N}} + \dfrac{dt}{2} \mat{A}^{j-1} \right]$
    \Comment{$\mat{I_{N}}$ is the $N \times N$ identity matrix.}
        \State $\mat{W} \gets \mat{W} \left[ \mat{dX}^j \right]^T $
    \State $\mat{C_P}^{K,J} \gets \mat{W}\ \mat{V}^j  $
    \State $\mat{C_P}^{J,K} \gets \left[  \mat{C_P}^{K,J} \right]^T  $
    \EndIf
\If{$j == 0$}
    \State $\mat{W} \gets \mat{\Lambda}^j \left[ \mat{V}^j \right]^T$
\ElsIf{$j < Nt$}     \State $\mat{W} \gets \text{vstack}\left(\mat{W}, \mat{\Lambda}^j \left[ \mat{V}^j \right]^T\right)$ \Comment{Vertically stack matrices.}
\EndIf
\EndFor
\Ensure $\mat{C_P}$, $\mat{V}^j$ $\forall j \in [0, N_t]$
\end{algorithmic}
\end{algorithm}

In the present work, $m$ is selected iteratively starting from a value of $m=1$ and increasing by 50\% at each iteration (rounding up).
At each iteration, an estimate of the total variance in the range at a specified DOF $\sigma_r^2$ is calculated based on up to 500 eigenpairs (denoted $N_{eig}$) of $\mat{C_P}$.
At the specified DOF, the difference between the maximum and minimum displacements across all times in the $i$th eigenvector of the DPSC is denoted as $r_i$. Then, the estimate of the total variance of the range of the DOF is estimated as
\begin{equation} \label{eq:estimated_variance}
    \sigma_r^2 = \sum_i^{N_{eig}} \lambda_{P,i} r_i^2,
\end{equation}
where $\lambda_{P,i}$ is the $i$th eigenvalue.
Note, that this assumes that all eigenpairs of $\mat{C_P}$ add independently to the total range experienced at a DOF. This is not generally true but is a useful approach for rapidly checking convergence of $\mat{C_P}$.
The value of $m$ is considered converged when
\begin{equation}
    \left| \dfrac{\sigma_r^2(m_{current}) - \sigma_r^2(m_{previous})}{\sigma_r^2(m_{current})}  \right| < 0.001.
\end{equation}

\FloatBarrier

\subsection{Sampling Direct Period-Spanning Covariance} \label{sec:sampling_DPSC}

The DPSC $\mat{C_P}$ can be used to directly sample the space of possible vibration responses over a full period. The eigenvectors of $\mat{C_P}$ correspond to independent perturbations from the mean (deterministic) response with variance equal to the corresponding eigenvalues. In general, using a subset of the eigenpairs with the highest eigenvalues may be sufficient for the statistical model.
Here, DPSC is sampled via Monte Carlo (DPSC-MC), assuming a normal distribution associated with each eigenpair. Latin hypercube sampling is employed.

\Cref{alg:dpsc_mc} outlines the sampling procedure for $N_s$ samples using $N_{eig}$ eigenpairs, the deterministic response at time $j$ of $\vect{x}^j$, a vector $\vect{R}$ defining the DOF of interest, and the reduced-order DPSC results from \Cref{alg:cov_across_time_rom}.
The result of \Cref{alg:dpsc_mc} is the matrix $\mat{x_{RS}}$ that has columns representing randomly sampled trajectories for the DOF of interest (specified by $\vect{R})$ and rows corresponding to different time instants.
Using $\mat{x_{RS}}$, the mean amplitude and range of amplitudes can easily be calculated for each sample and then used to construct a covariance between the mean cycle amplitude and the amplitude range over a cycle.

In the subsequent analysis, the number of eigenvalues used for MC, $N_{eig}$, is set to 500 for the large DPSC-MC computations. For the small DPSC-MC computation, $N_{eig}$ is dynamically set to cover at least an estimated 99.9\% of the standard deviation of the range of the response across a single cycle. This computation is done assuming that eigenpairs from $\mat{C_P}$ are distributed normally and independently contribute to the variance of the range (see \eqref{eq:estimated_variance}).
The independence assumption is conservative for estimating the variance of the cycle range because, generally, each eigenpair will not perfectly align the times of the maximum and minimum amplitude values to add to the cycle range variance. 
However, the independence assumption allows for a fast calculation to estimate the number of eigenpairs to use for the full DPSC-MC computation, which avoids the independence assumption from MC sampling.
The large DPSC-MC uses 1e6 samples, and the small DPSC-MC uses 1e3 samples.

\begin{algorithm}[h!]
\caption{DPSC-MC algorithm that samples stochastic trajectories using the DPSC. \indicesNotation}\label{alg:dpsc_mc}
\begin{algorithmic}[1]
\Require $N_t$, $m$, $N_{eig}$, $N_s$, $\mat{C_P}$, $\vect{x}^j$, $\mat{V}^j$ $\forall j \in [0, N_t]$, $\vect{R}$
\State $\vect{\lambda_P}, \mat{V_P} \gets \text{eigs}\left(\mat{C_P}, N_{eig}, \text{`largest'}\right)$ \Comment{Calculate $N_{eig}$ eigenpairs with largest eigenvalues.}
\State $K \gets 0$ to $N_{eig} - 1$ \Comment{Indices for all eigenvectors.}
\State $\mat{V_{PR}} \gets$ zeros($N_t + 1$, $N_{eig}$) \Comment{Eigenvectors extracted to DOF of interest.}
\State $\vect{x_R} \gets$ zeros($N_t+1$) \Comment{Time-varying mean response at DOF of interest.}
\For{$j=0$ to $N_t$}     \State $J \gets j m $ to $ (j+1) m - 1 $ \Comment{$J$ refers to the $j$th matrix block of size $m$.}
    \State $\mat{V_{PR}}^{j,K} \gets \vect{R}^T \left[\mat{V}^j\right]^T \mat{V_P}^{J,K}$
    \State $\vect{x_R}^j \gets \vect{R}^T \vect{x}^j$ 
\EndFor
\State $\mat{z} \gets $ standardNormal($N_{eig}$, $N_s$)
\State $\mat{\Sigma} \gets \text{diag}\left(\sqrt{\vect{\lambda_P}}\right) $ \Comment{Element-wise square root.}
\State $\mat{x_{RS}} \gets \vect{x_R} + \mat{V_{PR}}\ \mat{\Sigma} \  \mat{z} $ \Comment{Broadcast $\vect{x_R}$ to construct sampled trajectories of DOF.}
\Ensure $\mat{x_{RS}}$
\end{algorithmic}
\end{algorithm}

Given sample trajectories, the present work analyzes the time-averaged mean value of the cycle and the range of the response amplitudes of a cycle. For fatigue analysis, it is common to do rainflow counting of fatigue cycles and some form of racetrack filtering where only cycles that are large at a dominant frequency are used in the fatigue calculation \cite{veers_variable_amplitude_1989, hayman_mlife_2012}. The smaller cycles still affect the range of amplitudes experienced per large cycle, but not the number of cycles. For the present work, the larger cycles correspond to the deterministic solution, and the stochastic excitation and response create small deviations that influence the cycle range and mean.
Statistics for the time-averaged cycle mean and cycle range are calculated from the samples of these quantities derived from the MC trajectories. 
Additional methods for counting fatigue cycles are detailed in a standard \cite{astm_fatigue_counting} and could be applied within the current framework using full cycles of the response sampled from DPSC-MC.
Frequency domain methods also exist for analyzing fatigue damage but may not as easily incorporate nonlinear or time-periodic influences considered here \cite{dirlik_dirlik_2021}.

While time-averaged mean and range per cycle are used to analyze the systems here, full stochastic trajectories for individual cycles are produced as an intermediate step, so other cycle properties could easily be used instead within the DPSC-MC framework. For instance, the maximum displacement for any time on a cycle could be used to check the probability of exceeding a displacement limit per cycle. This cannot be directly calculated from the DPC because that does not provide a direct way to determine the time at which the maximum amplitude occurs. The maximum amplitude may not always be exceeded at the same time, and DPC does not capture correlations between time instants.

The benefit of the present approach is that applying Monte Carlo to sample stochastic trajectories of full cycles at a time is much faster than doing nonlinear stochastic time integration for an equivalent number of cycles. Furthermore, the time integration must generally be serialized while stochastic trajectories can be independently sampled.
Preliminary work considered polynomial chaos expansion rather than random sampling to evaluate statistics, but polynomial chaos expansion showed poor convergence. This is possibly due to the non-smooth (piecewise planar) nature of the maximum amplitude with respect to the many different eigenpairs of $\mat{C_P}$.

\FloatBarrier

\section{Extension of 2-DOF Example} \label{sec:2dof_example}

First, the DPSC is demonstrated for the second example from \cite{ahsan_adjoint_based_2025}. The 2-DOF system has time derivative 
\begin{equation}
    \vect{\dot{x}} = 2 \pi \begin{bmatrix}
        x_2 \\ -2 x_2 - x_1 + 2 \cos(2 \pi t)
    \end{bmatrix}
\end{equation}
with stochastic forcing applied as
\begin{equation}
    \mat{F} = \begin{bmatrix}
        0
        \\ x_1
    \end{bmatrix}.
\end{equation}

\begin{figure}[h!]
\centering
        \includegraphics[scale=1,width=0.9\linewidth]{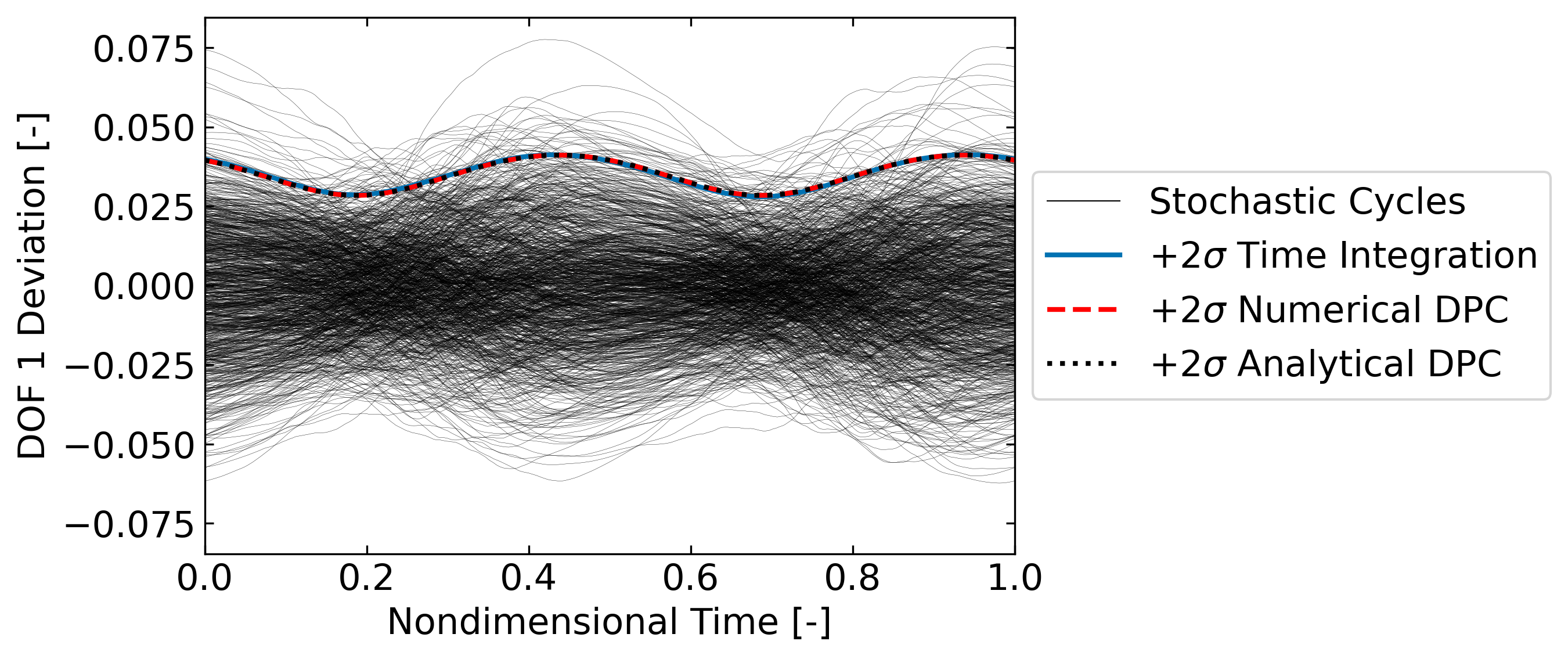} 
    \caption{DOF 1 time-varying deviation from deterministic solution. Only the last 1000 cycles of the time integration response are shown.}
    \label{fig:ahsan2_std}
\end{figure}

\begin{figure}[h!]
\centering
        \includegraphics[scale=1,width=0.7\linewidth]{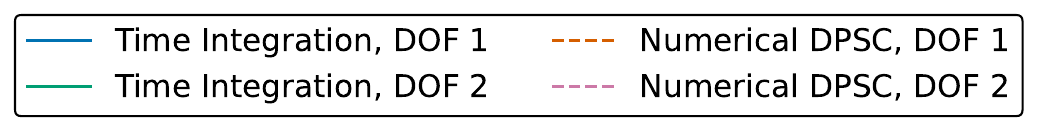} 
    \\
	\begin{subfigure}[b]{0.49\linewidth } 
		\centering
		\includegraphics[width=\linewidth]{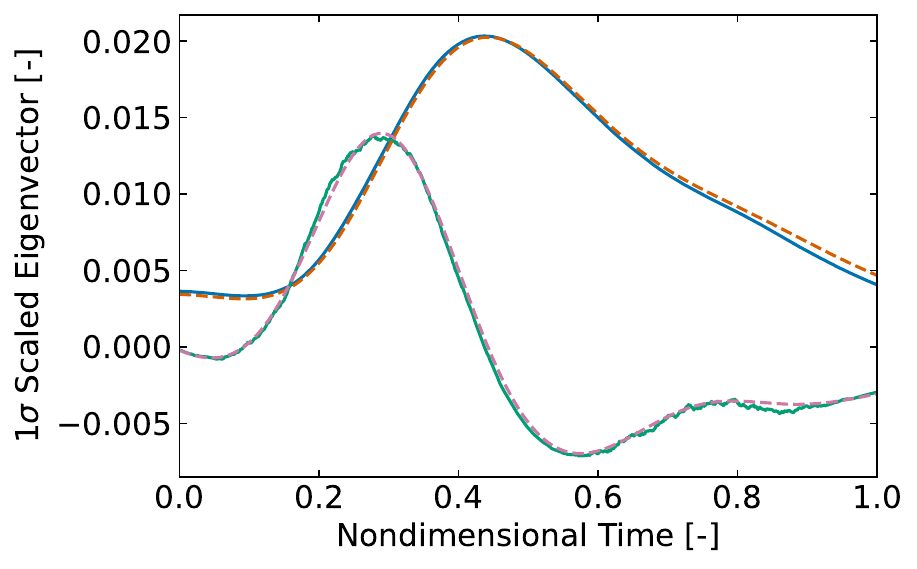}
		\caption{}
        \label{fig:ahsan2_period_cov_sub_a}
	\end{subfigure}	\hfill
	\begin{subfigure}[b]{0.49\linewidth } 
		\centering
		\includegraphics[width=\linewidth]{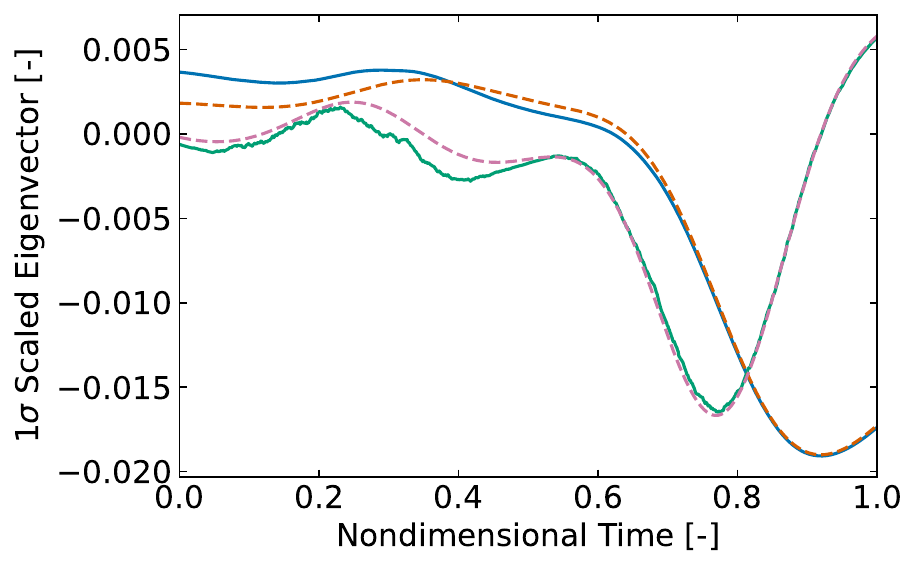}
		\caption{}
	\end{subfigure}    \\
	\begin{subfigure}[b]{0.49\linewidth } 
		\centering
		\includegraphics[width=\linewidth]{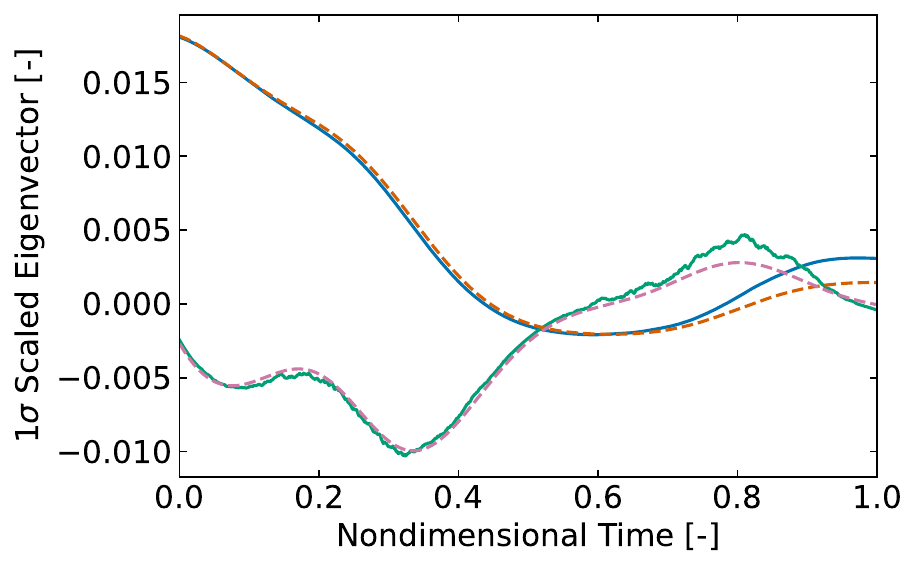}
		\caption{}
	\end{subfigure}	\hfill
	\begin{subfigure}[b]{0.49\linewidth } 
		\centering
		\includegraphics[width=\linewidth]{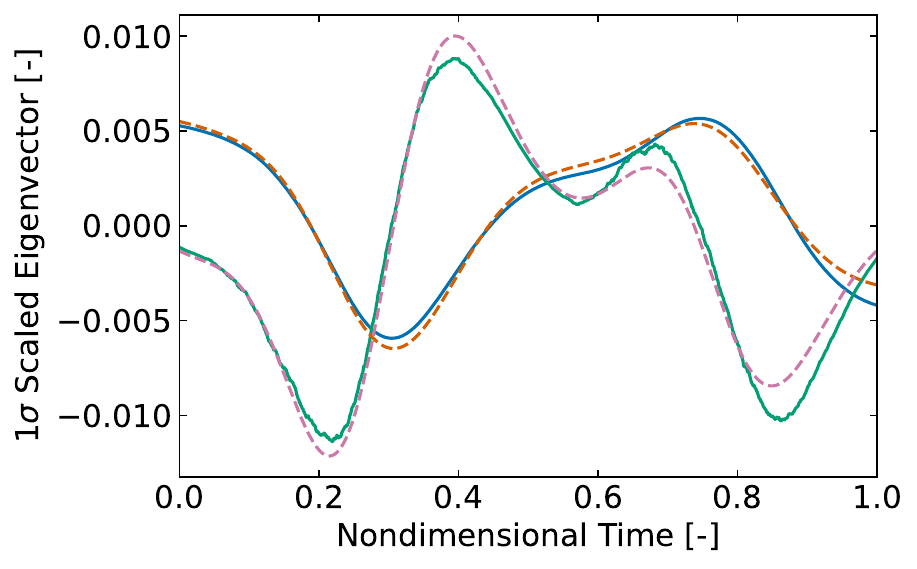}
		\caption{}
	\end{subfigure}    \caption{Time representations of period-spanning covariance eigenvectors with indices (a) 1, (b) 2, (c) 3, (d) 4.}
    \label{fig:ahsan2_period_cov}
\end{figure}

\Cref{fig:ahsan2_std} shows that the DPC from the modified algorithm is consistent with the time integration (4000 cycles, 2$^\text{14}$=16,384 steps per cycle) standard deviation. Both time-varying standard deviations are also consistent with the analytical solution presented in \cite{ahsan_adjoint_based_2025}.
The extension to the DPSC in this work is shown in \Cref{fig:ahsan2_period_cov}, where the first four eigenvectors are compared between the DPSC computation and calculations from the time integration. 
These eigenvectors represent the largest common perturbations caused by the white noise excitation that are correlated across a full vibration cycle.
These eigenvectors are not periodic because system responses to excitation decay based on the system properties, so they are not expected to be periodic at the external forcing frequency. 
Rather, the eigenvectors of the DPSC provide a convenient space of independent vectors to be used in sampling probable stochastic responses.

The numerical DPSC solution in \Cref{fig:ahsan2_period_cov} uses the analytical solution for the deterministic problem and then numerically calculates all of the covariance solution steps.
For the first eigenvector (see \Cref{fig:ahsan2_period_cov_sub_a}), the match is extremely good. For the other eigenvectors, there are some errors in the match between the DPSC and time integration solutions. However, these are attributed to insufficient convergence of the time integration solution.
Given this verification of the DPSC, two larger systems are considered next using the DPSC-MC approach to directly sample the stochastic cycles.

\FloatBarrier

\section{Aeroelastic System} \label{sec:aero_full}

To test the methodology for larger systems approaching engineering relevance, this section considers an aeroelastic system consisting of a beam and unsteady aerodynamics. \Cref{sec:aero_describe} describes the system.
Results are then presented for the response in \Cref{sec:aero_resp_results} with computation times in \Cref{sec:aero_comp_time}.
Additional insights into the DPSC are presented in \Cref{sec:aero_period_spanning}.
Convergence of the stochastic simulations and conditioning of the DPC method are discussed in Sections \ref{sec:aero_convergence} and \ref{sec:aero_conditioning}, respectively.

\subsection{System Description} \label{sec:aero_describe}

The aeroelastic system consists of a linear structural model with constant properties and a constant airfoil with unsteady aerodynamics.
The intent of the system is not to exactly model any specific physical system but rather to demonstrate the potential of the methods to scale to engineering-relevant aeroelastic systems.
The system is represented with a variable number of finite elements to investigate the scaling of computational costs. 
The structural model consists of a cantilevered beam and uses Euler-Bernoulli linear finite beam elements with cubic shape functions in the flap (perpendicular to airfoil chord) and edge (parallel to airfoil chord) directions. Linear finite elements are used for torsion. The beam is 100~m long, has a mass length density of 550~kg/m, and a rotational inertia length density of 6875~kgm. Bending and torsional stiffnesses are set to have analytical first natural frequencies in the flap-bending, edge-bending, and torsion directions of 0.5 Hz, 0.7 Hz, and 4.0 Hz, respectively.
The structural model uses mass and stiffness proportional damping coefficients of 0.04~1/s and 0.001~s, respectively.

The unsteady aerodynamics model from \cite{branlard_dynamic_2022} is applied at 4 Gauss quadrature points per finite element and describes the history-dependent evolution of the aerodynamic forces at each point.
Aerodynamic parameters for a FFA-W3-211 airfoil are presented in \Cref{tab:aero_pars}, and full details can be found at \cite{porterStochasticCodeData2026}.
Parameters are taken from \cite{branlard_dynamic_2022} and associated code \cite{branlardWelibZenodo2022}.
The unsteady aerodynamics model has four first-order states per evaluation point that describe the hysteretic evolution of the aerodynamic forces and includes a Wagner's memory function \cite{branlard_dynamic_2022}.
The exact details of the aerodynamic model are not of primary interest here; rather, this work demonstrates that the present methods work for aerodynamic models relevant to engineering research.
The present implementation \cite{porterStochasticCodeData2026} can be run to reproduce some plots from \cite{branlard_dynamic_2022} to verify the current unsteady aerodynamic model. An example plot of the hysteretic aerodynamic forces is included next in \Cref{sec:aero_resp_results}.

\FloatBarrier

\begin{table}[h!]
\centering
\caption{Unsteady aerodynamics parameters based on \cite{branlard_dynamic_2022,branlardWelibZenodo2022}.} \label{tab:aero_pars}
\begin{tabular}{clc}
\hline\hline
Symbol & Description    & Value                         \\ \hline
$\rho_{air}$ & Density of air & 1.225 kg/m$^3$ \\ 
$c$ & Chord          & 3.0 m                         \\
$x_{pitch}$ & Pitch axis as fraction of chord       & 0.5                           \\
$\alpha_0$ & Zero lift angle of attack         & -2.99 degrees                 \\
$C_{l,\alpha}$ & Lift slope at zero lift angle of attack            & 7.376                         \\
$A_1$ & Parameter of Wagner's memory function            & 0.165                         \\
$A_2$ & Parameter of Wagner's memory function            & 0.335                         \\
$b_1$ & Parameter of Wagner's memory function           & 0.0445                        \\
$b_2$ & Parameter of Wagner's memory function            & 0.3                           \\
$T_{f,0}$ & Coefficient for trailing-edge stall time constant          & 6.0                           \\
$T_{p,0}$ & Coefficient for boundary layer pressure gradient           & 1.5                           \\
\hline\hline
\end{tabular}
\end{table}

External cosine excitation is applied with tip loads of 540~kN in the edge direction and 100 kNm in the torsional direction at a forcing frequency of 0.127 Hz (0.8 rad/s). 
The system is also subject to an inflow velocity with a mean of 40 m/s and harmonic oscillations of $\pm$26 m/s in phase with the external forcing. 
Harmonic oscillations in the inflow can be caused by the rotation of a structure subject to flow velocity that is not parallel to the axis of rotation (e.g., for a helicopter rotor).
The nominal angle of attack of the airfoil is 15 degrees, but the instantaneous angle of attack is calculated in the unsteady aerodynamics model based on the local velocities of the structure.
Steady-state solutions are found via time integration and a nonintrusive HBM wrapper as described in Appendix~\ref{sec:appendix_methods}.

\FloatBarrier

\subsection{Aeroelastic Response Results} \label{sec:aero_resp_results}

\begin{figure}[h!]
\centering
        \includegraphics[scale=1,width=\linewidth]{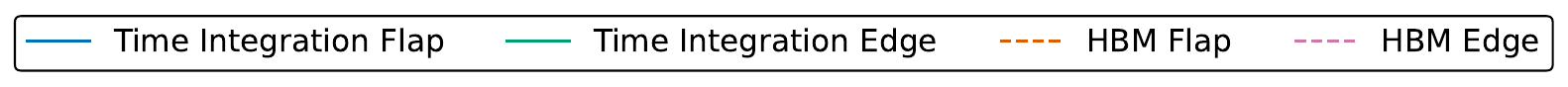} 
	\begin{subfigure}[b]{0.49\linewidth } 
		\centering
		\includegraphics[width=\linewidth]{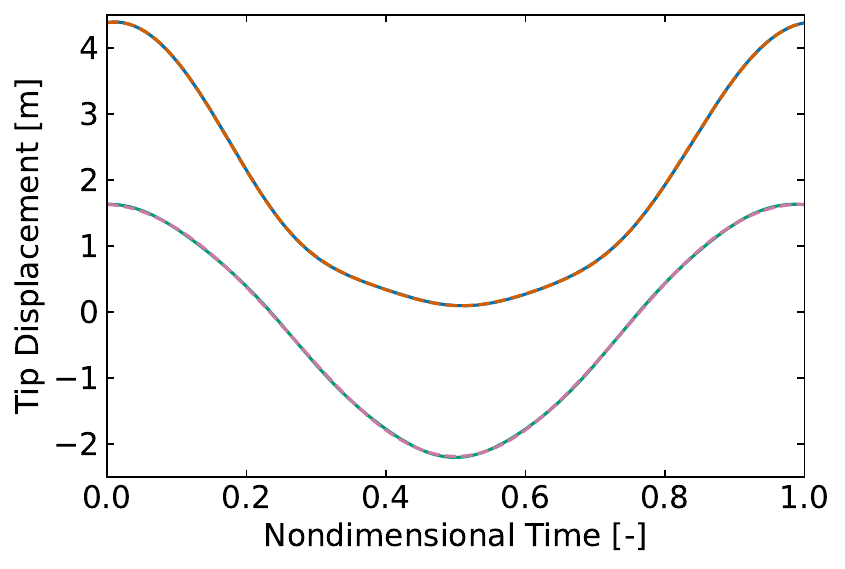}
		\caption{}
	\end{subfigure}	\hfill
	\begin{subfigure}[b]{0.49\linewidth } 
		\centering
		\includegraphics[width=\linewidth]{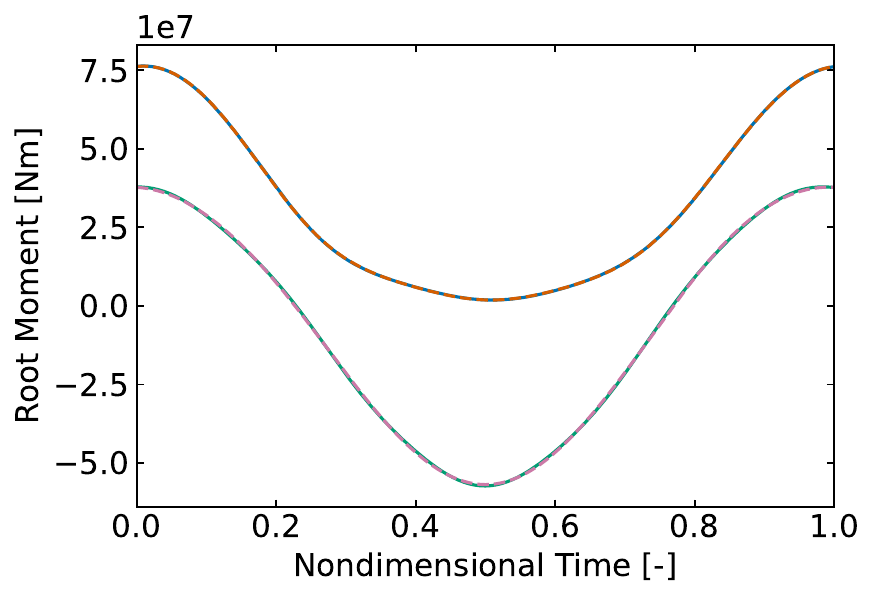}
		\caption{}
	\end{subfigure}\caption{Steady-state solutions for 16 elements for flap and edge (a) tip displacements and (b) root bending moment.}
\label{fig:ba_steady}
\end{figure}

\begin{figure}[h!]
\centering
    	\begin{subfigure}[b]{0.49\linewidth } 
		\centering
		\includegraphics[width=\linewidth]{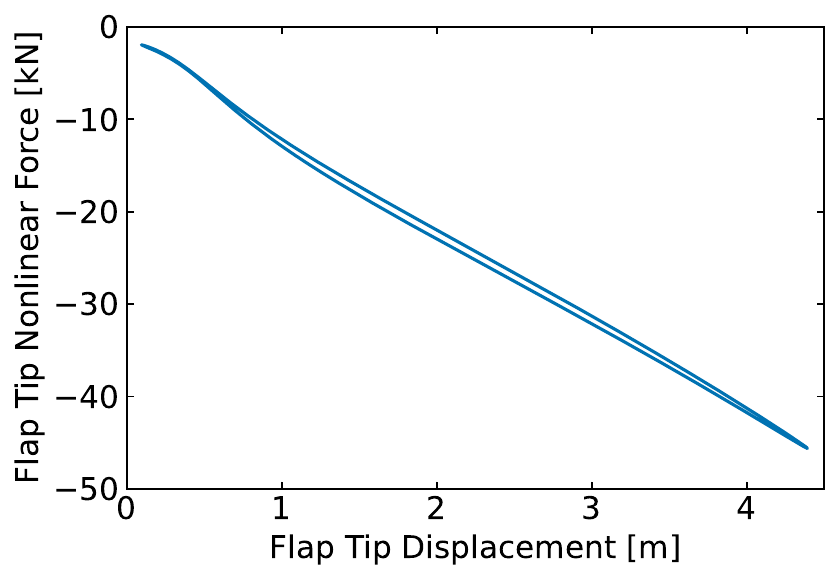}
		\caption{}
	\end{subfigure}	\hfill
	\begin{subfigure}[b]{0.49\linewidth } 
		\centering
		\includegraphics[width=\linewidth]{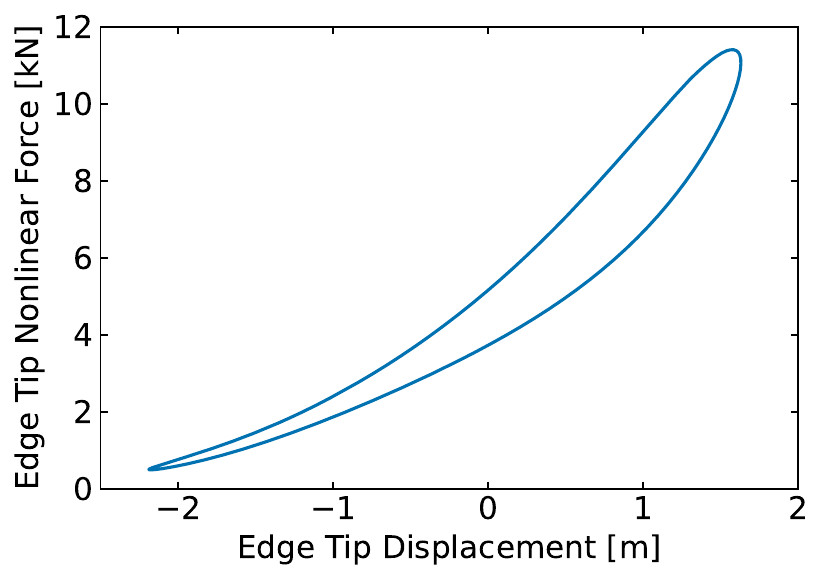}
		\caption{}
	\end{subfigure}    \caption{Steady-state aerodynamic forces acting on the tip node for 16 elements (a) flap direction and (b) edge direction as calculated from the HBM solution.}
    \label{fig:aero_forces}
\end{figure}

This section reviews the results for the nominal model using 16 finite elements.
\Cref{fig:ba_steady} shows the steady-state responses as calculated from both time integration and HBM.
For both the tip displacements and the root moments the two approaches yield essentially identical responses, indicating that the HBM solution is converged when using the zeroth and first five harmonics and 128 AFT steps. The computational benefits of the HBM solution are highlighted later in \Cref{sec:aero_comp_time}. The flap responses also show some clear distortions from a simple single harmonic cosine response as a result of the nonlinear effects of the varying inflow and the unsteady aerodynamics.
Additionally, the integrated effect of the aerodynamic forces on the tip DOF are highlighted in \Cref{fig:aero_forces} and show clear hysteresis. The evolution of the nonlinear forces is due to a combination of the time-varying inflow rate (which is highest at the time zero corresponding the largest positive displacements in both directions) and the hysteretic evolution of the nonlinear forces.

\FloatBarrier

\begin{figure}[h!]
\centering
        \includegraphics[scale=1,width=0.8\linewidth]{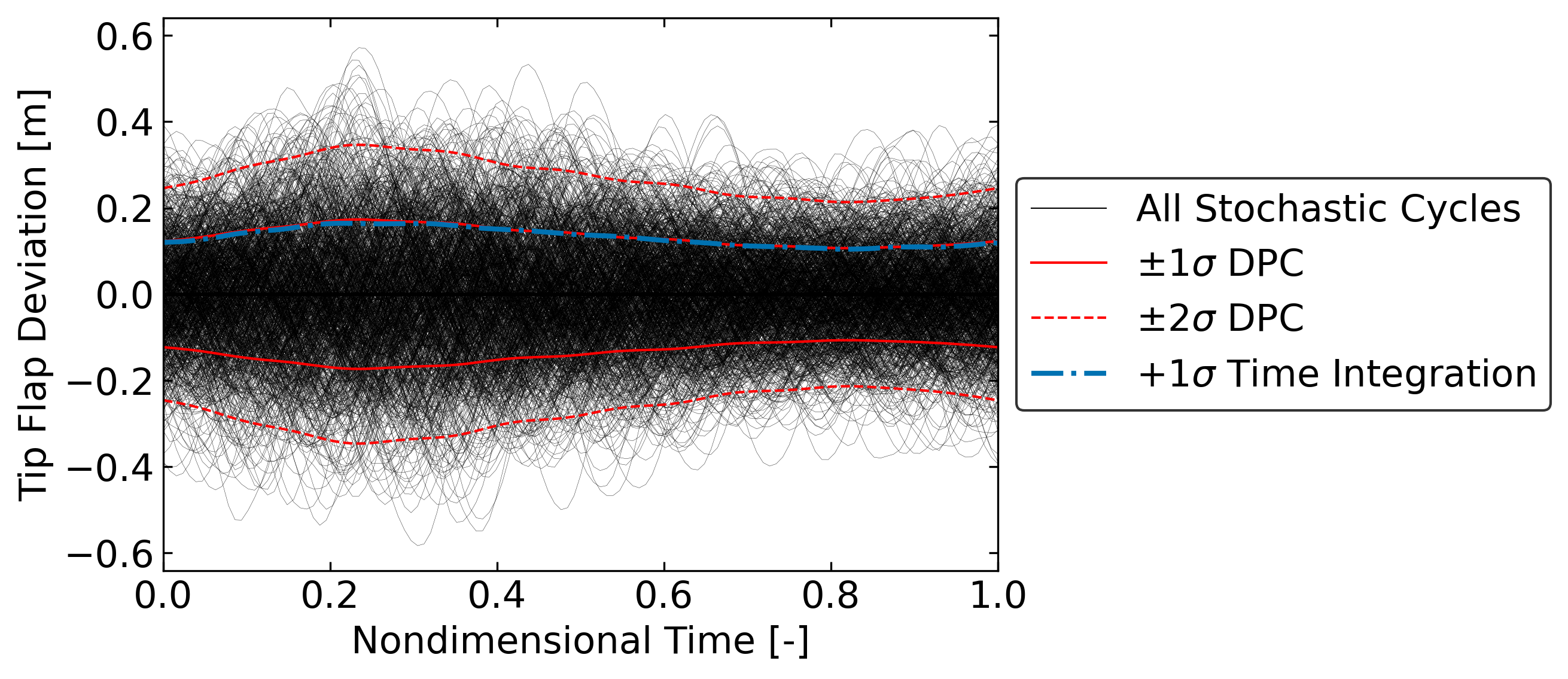} 
\caption{Tip flap motion comparing stochastic trajectories, standard deviation from time series, and standard deviation from direct adjoint calculation.}
\label{fig:ba_tip_flap_std}
\end{figure}

Building on the steady-state solution, the effect of white noise excitation on the response is shown by the deviations from the steady-state solution in \Cref{fig:ba_tip_flap_std}. The white noise is caused by inflow variations with $\sigma = 0.4$ m/s.
Here, the DPC method exactly captures the standard deviation of the 1000 cycles of time integration of the stochastic responses. Time integration is conducted with exactly 128 time steps per cycle, so standard deviations can be calculated from 1000 points for each sampled phase angle of the external forcing in \Cref{fig:ba_tip_flap_std}. The DPC calculation also discretizes the period into 128 steps.
A time-varying standard deviation in the response is observed rather than a constant standard deviation because of how the inflow perturbations interact with the time-varying inflow and the periodically varying position and velocity of the beam. Thus, methods for calculating a constant standard deviation of the response will miss these time-varying influences.
However, a time-varying standard deviation is not necessarily directly applicable to engineering design considerations.

\begin{figure}[h!]
\centering
        \includegraphics[scale=1,width=0.7\linewidth]{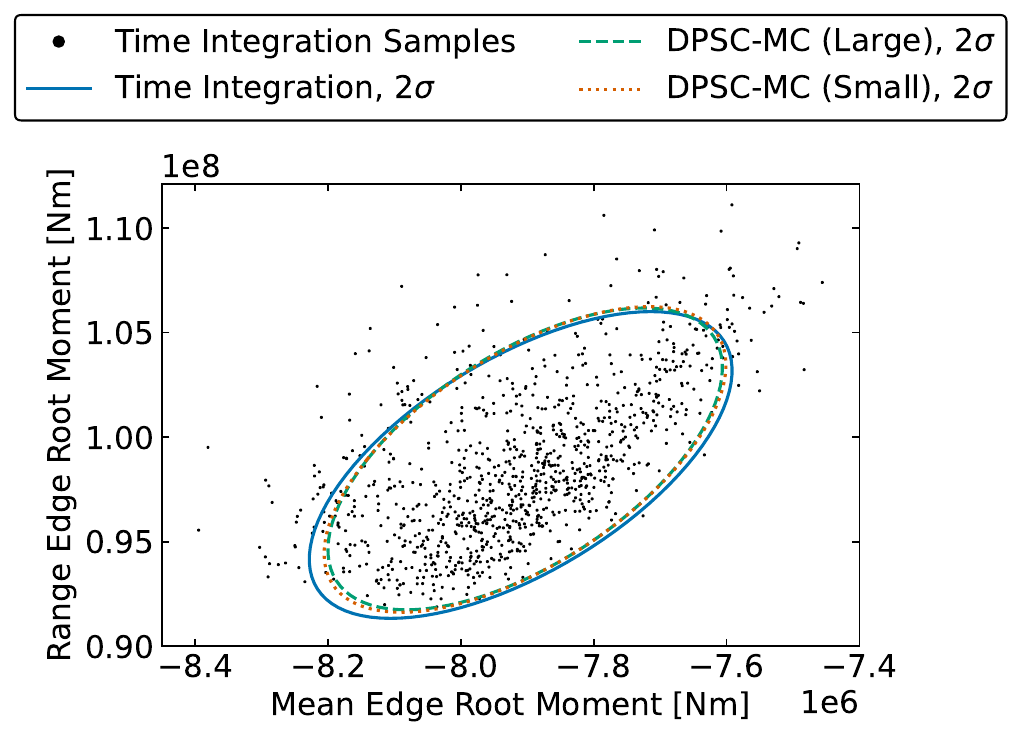} 
    \caption{Covariance of mean and range over a cycle of the root bending moment in the edge direction.}
    \label{fig:ba_edge_moment_std}
\end{figure}

Rather, to estimate fatigue life, the mean amplitude and range of amplitudes for each cycle need to be sampled for the stochastic response. 
\Cref{fig:ba_edge_moment_std} shows the covariance between these two response variables for the root bending moment in the edge direction. The present calculations from DPSC-MC accurately capture the covariance that is estimated from the time integration. Furthermore, the DPSC-MC method allows for direct sampling in the 2D space shown in \Cref{fig:ba_edge_moment_std}, so it can provide further details beyond just the first two moments of the distribution of fatigue responses (see \Cref{sec:aero_convergence}).
In \Cref{fig:ba_edge_moment_std}, the small and large MC covariances are consistent, indicating that the small MC covariance uses sufficient samples (at lower computational cost). The slight differences between the DSPC-MC and time integration covariances are attributed to insufficient time integration cycles, as highlighted by the convergence investigations in \Cref{sec:aero_convergence}.
Next, the computational benefits of the present approach are highlighted.

\FloatBarrier
\subsection{Computation Time} \label{sec:aero_comp_time}
\newcommand{\hbmSpeedupMin}[0]{17$\times$}
\newcommand{\hbmSpeedupMax}[0]{26$\times$}
\newcommand{\stochSpeedupMin}[0]{67$\times$}
\newcommand{\stochSpeedupMax}[0]{142$\times$}
\newcommand{\stochShortSpeedupMin}[0]{292$\times$}
\newcommand{\stochShortSpeedupMax}[0]{1393$\times$}
\FloatBarrier

Computation time comparisons are made by running simulations on a dedicated compute node using high performance computing resources at the National Laboratory of the Rockies. All compared simulations are run during the same job submission on the same node.
Compute nodes have 104 CPU cores at 2.1 GHz and at least 256 GB of RAM. Some computations exploit parallelism, but the full computational power of the compute nodes is not necessary. Rather, a full compute node is just used to ensure that there are not outside influences on the computation time.

\begin{figure}[h!]
\centering
    	\begin{subfigure}[b]{0.49\linewidth } 
		\centering
		\includegraphics[width=\linewidth]{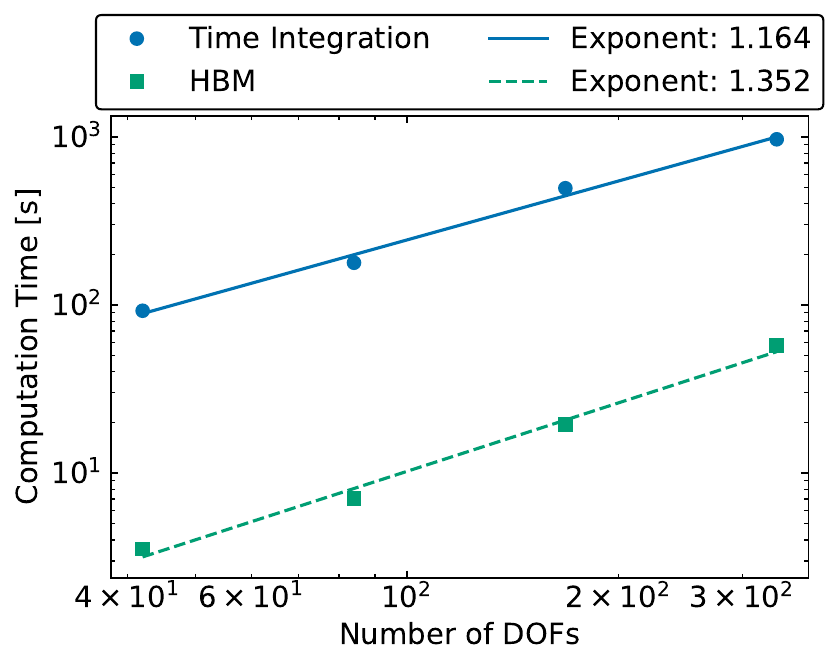}
		\caption{}
        \label{fig:ba_steady_time}
	\end{subfigure}	\hfill
	\begin{subfigure}[b]{0.49\linewidth } 
		\centering
		\includegraphics[width=\linewidth]{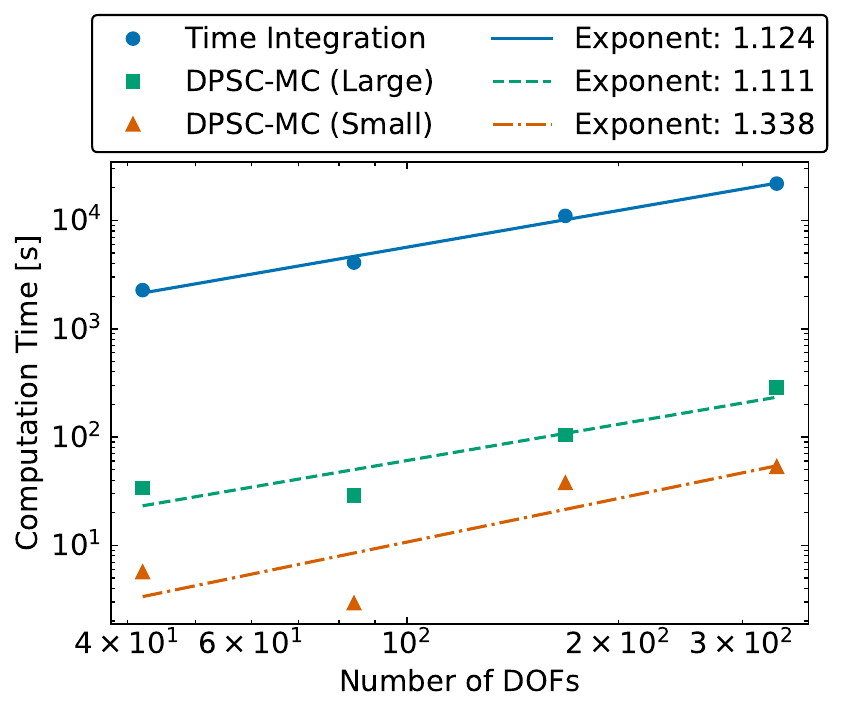}
		\caption{}
        \label{fig:ba_stochastic_time}
	\end{subfigure}\caption{Computation time for models with 2 to 16 elements for (a) steady-state solution and (b) stochastic calculations.}
\label{fig:ba_comp_time}
\end{figure}

The present methodology is significantly faster than stochastic time integration. 
For the steady-state responses, HBM is {\hbmSpeedupMin} to {\hbmSpeedupMax} faster than time integration of 40 cycles to approximate the steady-state solution. This speedup for HBM is expected compared to time integration \cite{krackHarmonicBalanceNonlinear2019}.
None of the time integration cases reached a relative error of less than 0.001 on the norm of the displacement states, velocity states, and first-order aerodynamic states, matching between the start and end of an integrated cycle of forcing.
As shown in \Cref{fig:ba_steady_time}, both time integration and HBM scale roughly linearly in computational cost with the number of DOFs. This is likely due to the linear increase in nonlinear force evaluations with the number of DOFs and relatively small matrix computations.

\FloatBarrier

Of more interest here, the DPSC-MC with the large number of MC samples and convergence checks is {\stochSpeedupMin} to {\stochSpeedupMax} faster than time integration of 1000 cycles. Furthermore, if convergence checks are not done, the DPSC-MC method is {\stochShortSpeedupMin} to {\stochShortSpeedupMax} faster than time integration while still giving converged results. Given that the 1000 cycles of time integration are not fully converged in \Cref{sec:aero_resp_results}, it is fair to compare to the converged small DPSC-MC results.
Here, the speedup by using DPSC-MC is equivalent to only needing to do 0.7--3.4 cycles of stochastic time integration. Thus, the present methodology represents a significant benefit over stochastic time integration.
For the stochastic simulations, the computational scaling remains similar for both methods and slightly over an exponent of 1. DPSC-MC uses larger matrices than time integration because the second-order states are converted to first-order states, but this represents a fixed coefficient difference on the computational cost of matrix operations and not a different exponent in computational scaling.
Importantly, a discrete Lyapunov solution algorithm is used to solve the matrix problem for $\mat{C_0}$ in DPC rather than a simpler approach of constructing an $N^2 \times N^2$ matrix, which would introduce poor scaling of computational costs and memory requirements.

\FloatBarrier
\subsection{Period-Spanning Covariance} \label{sec:aero_period_spanning}
\FloatBarrier

\begin{figure}[h!]
\centering
        \includegraphics[scale=1,width=0.9\linewidth]{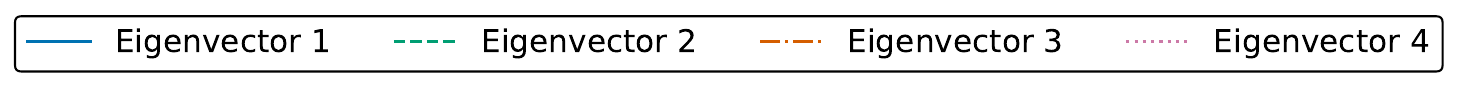} 
	\begin{subfigure}[b]{0.49\linewidth } 
		\centering
		\includegraphics[width=\linewidth]{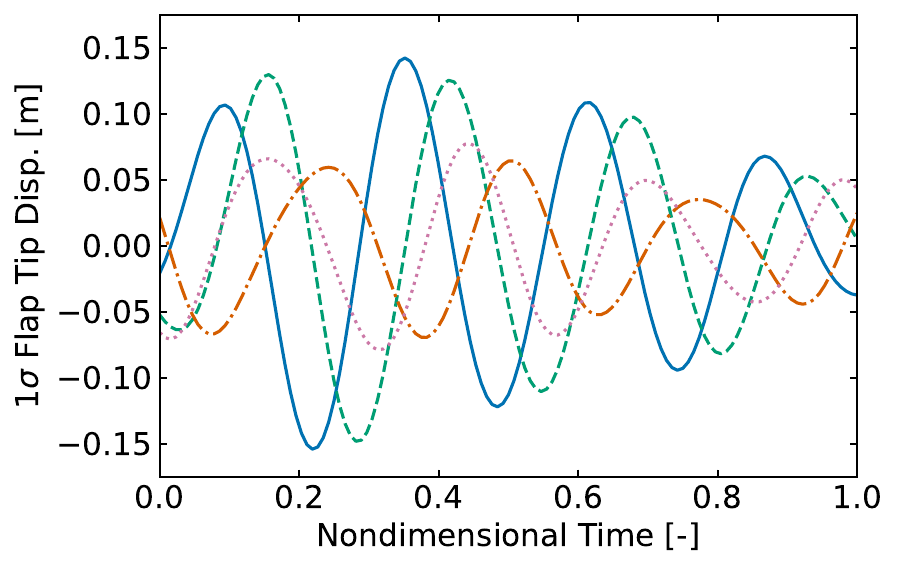}
		\caption{}
	\end{subfigure}	\hfill
	\begin{subfigure}[b]{0.49\linewidth } 
		\centering
		\includegraphics[width=\linewidth]{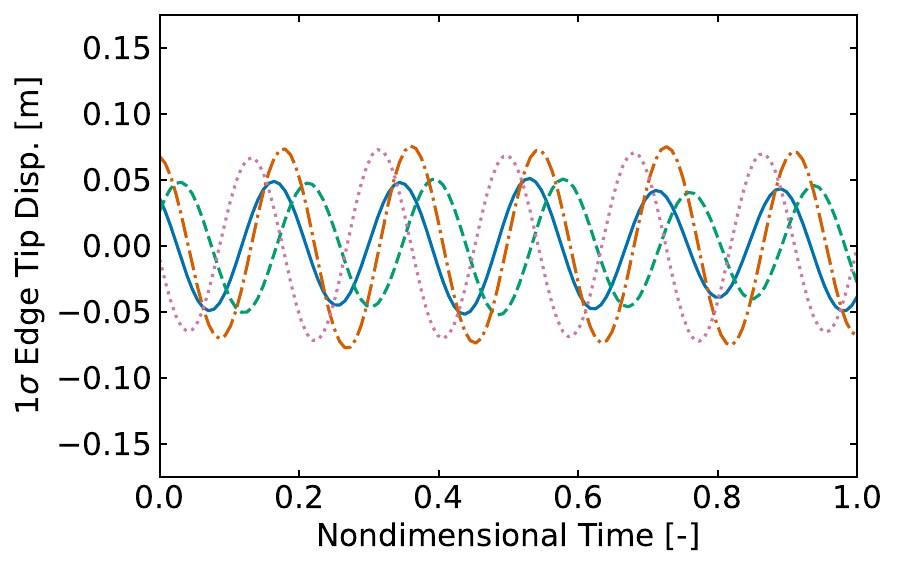}
		\caption{}
	\end{subfigure}\caption{Time representations of DPSC eigenvectors for tip displacements in (a) flap and (b) edge directions.}
\label{fig:ba_cov_time}
\end{figure}

This section briefly discusses the eigenvectors of the DPSC matrix to strengthen the conceptual understanding of the DPSC-MC method.
\Cref{fig:ba_cov_time} shows the time-varying tip displacements for the first four eigenvectors of the DPSC matrix. 
The flap tip displacements of the first four eigenvectors appear to oscillate at slightly below the first flap bending natural frequency. The edge tip displacements of the first four eigenvectors all appear to oscillate at near the first edge bending natural frequency.
For both cases, the responses are at different phases, so, collectively, the eigenvectors can represent probable responses while not individually representing single modes.

\FloatBarrier

\begin{figure}[h!]
\centering
        	\begin{subfigure}[b]{0.49\linewidth } 
		\centering
		\includegraphics[width=\linewidth]{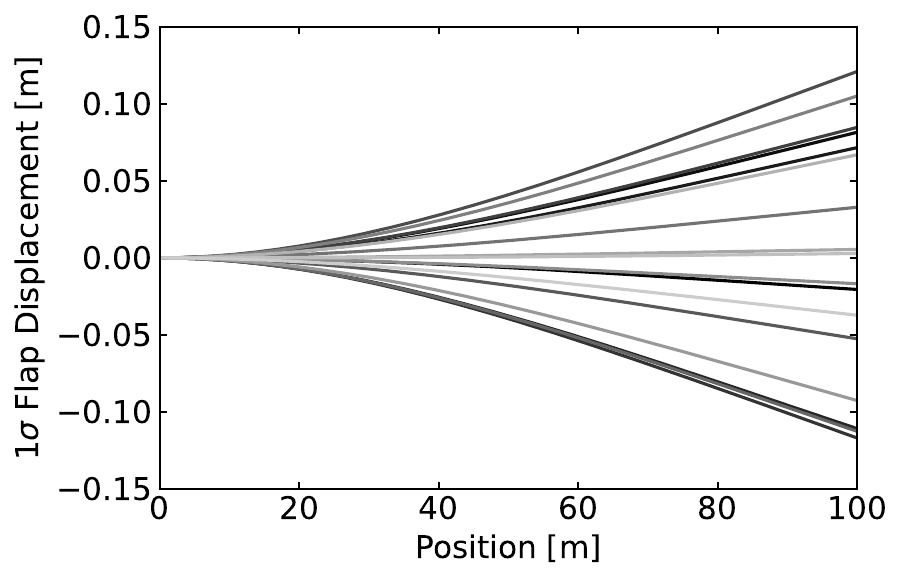}
		\caption{}
	\end{subfigure}	\hfill
	\begin{subfigure}[b]{0.49\linewidth } 
		\centering
		\includegraphics[width=\linewidth]{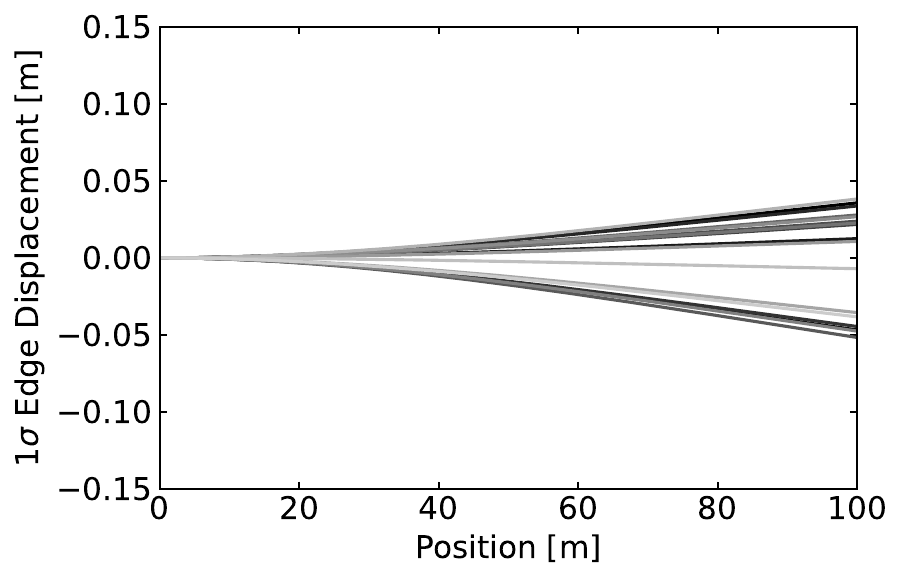}
		\caption{}
	\end{subfigure}\caption{Spatial representations of first eigenvector of the DPSC for (a) flap displacements and (b) edge displacements. Through the cycle, colors vary from black (initial time) to light gray (final time).}
\label{fig:ba_cov_space}
\end{figure}

These covariance eigenvectors also have spatial components. For the first four, the spatial distributions are similar to classic first bending mode shapes (see \Cref{fig:ba_cov_space}). Beyond the fourth eigenvector, the shapes are not easily interpretable.
Because the eigenvectors include both time and spatial information, full stochastic cycles of the response with time and spatial information can be sampled with DPSC-MC and then reduced to any quantity of interest.

\FloatBarrier

\subsection{Convergence} \label{sec:aero_convergence}

\begin{figure}[h!]
\centering
        \includegraphics[scale=1,width=0.8\linewidth]{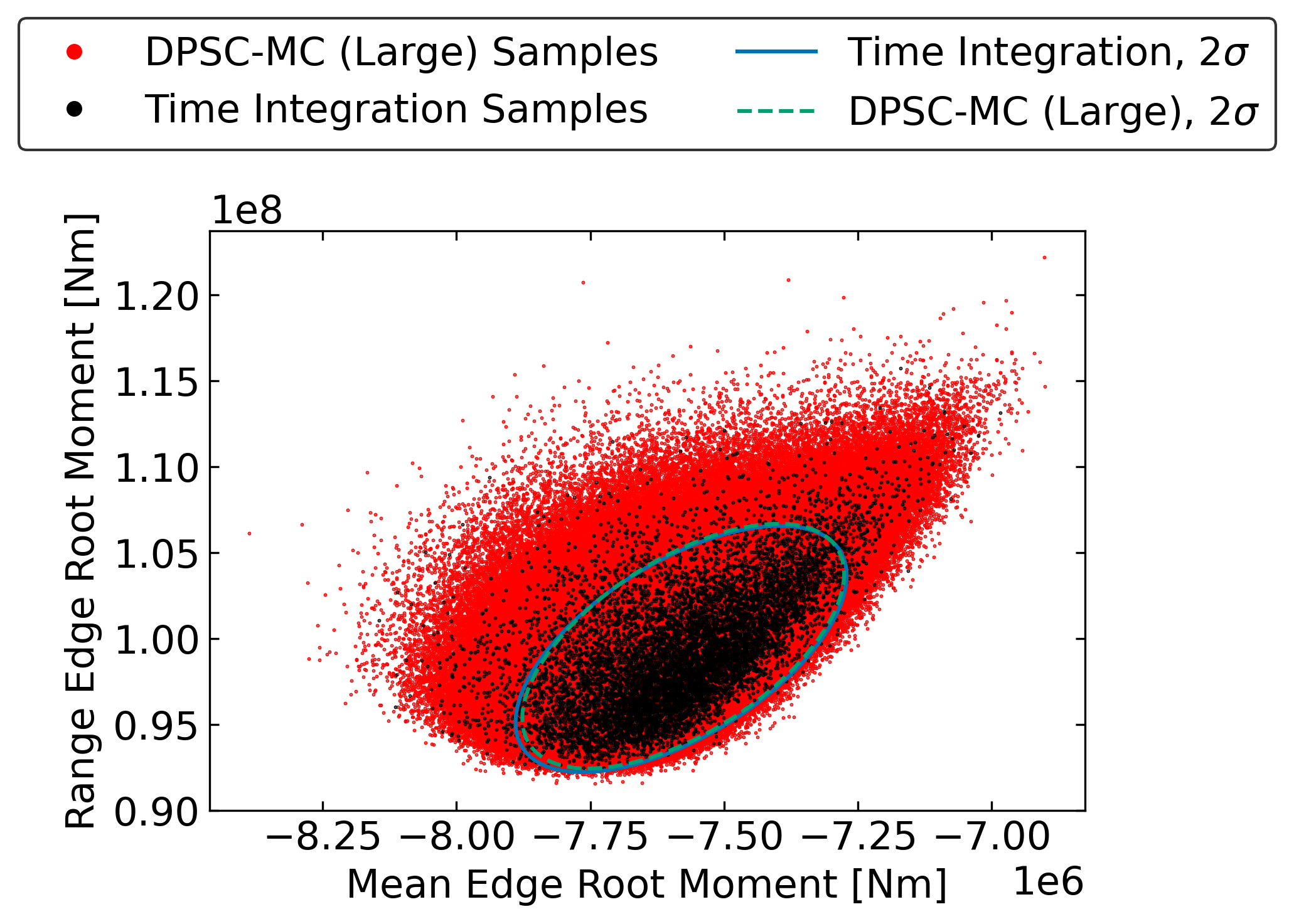} 
\caption{Joint distribution of cycle mean and range calculated from time integration and Monte Carlo based on the direct method for the edge direction root bending moment.}
\label{fig:ba_cov_dist_convg}
\end{figure}

\begin{figure}[h!]
\centering
        \includegraphics[scale=1,width=0.5\linewidth]{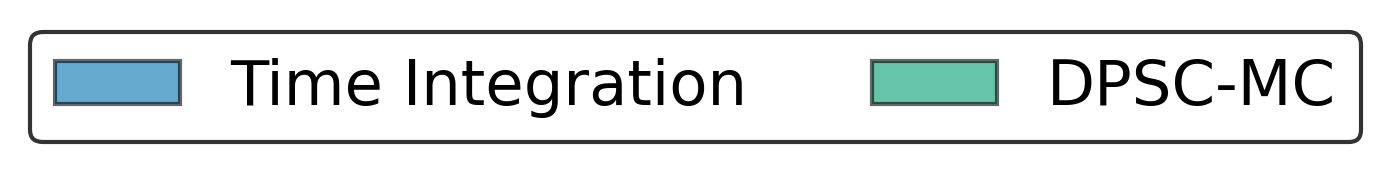} \\
	\begin{subfigure}[b]{0.49\linewidth } 
		\centering
		\includegraphics[width=\linewidth]{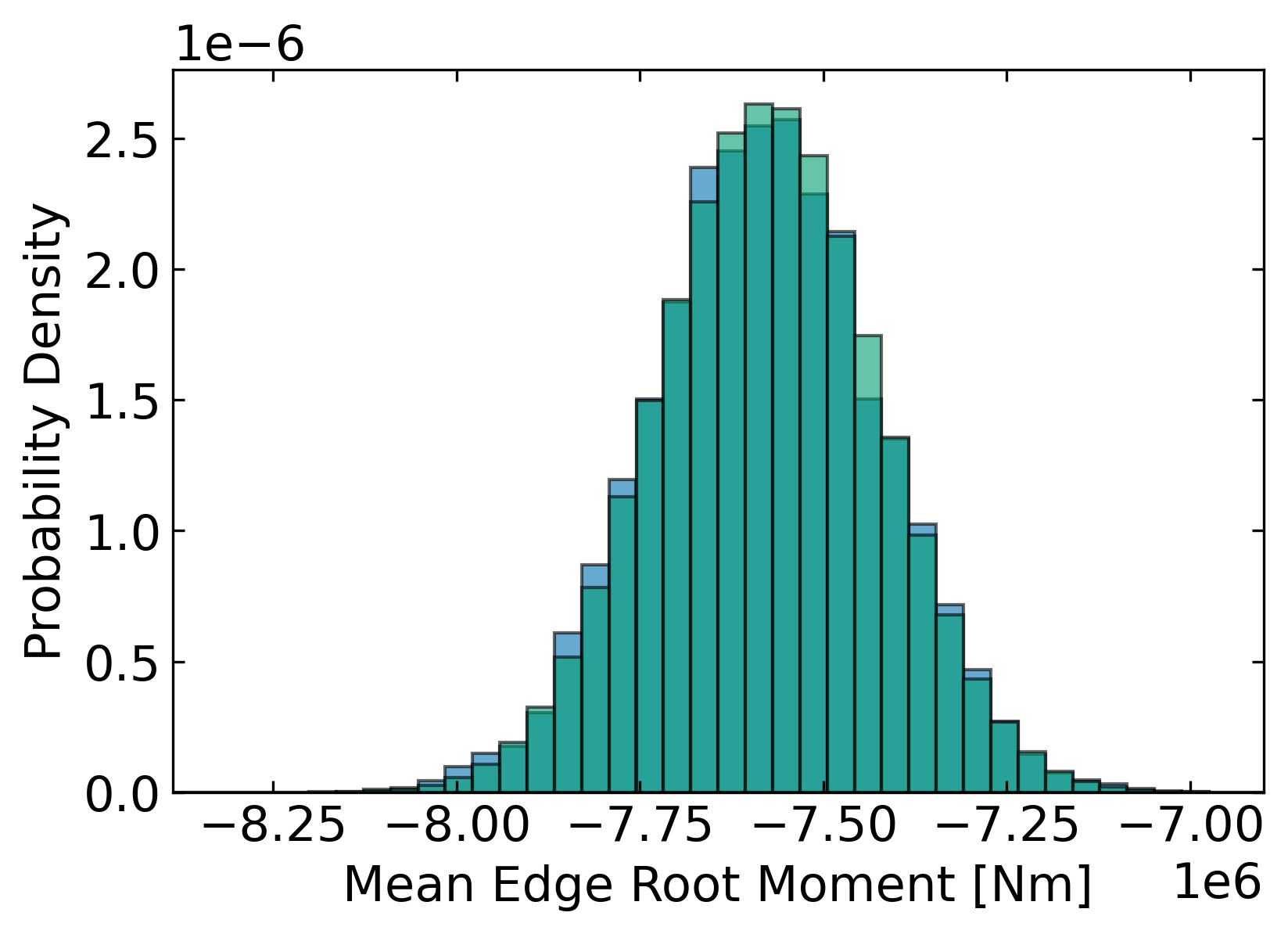}
		\caption{}
	\end{subfigure}	\hfill
	\begin{subfigure}[b]{0.49\linewidth } 
		\centering
		\includegraphics[width=\linewidth]{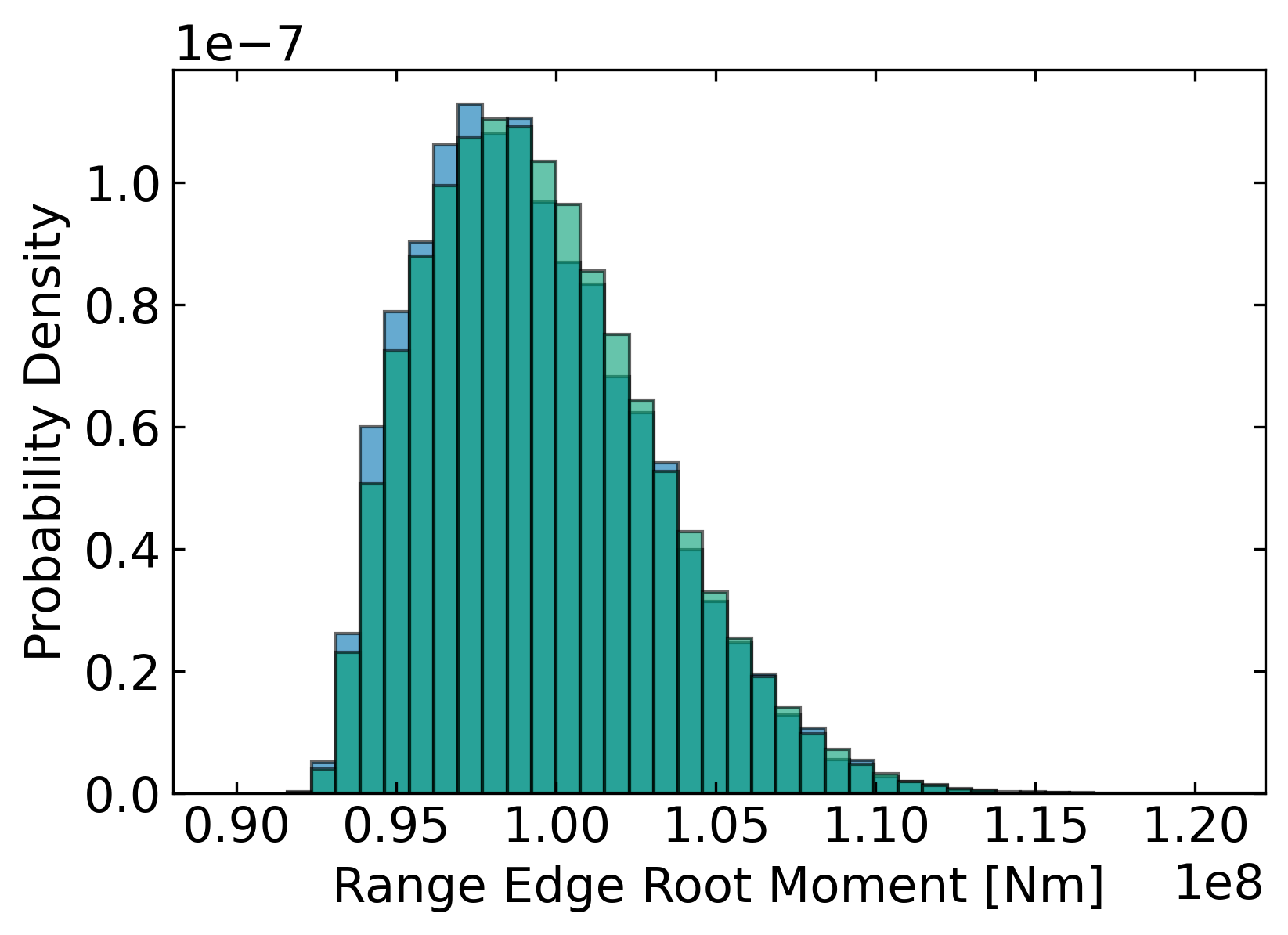}
		\caption{}
	\end{subfigure}\caption{Distributions of edge root moment (a) cycle mean and (b) cycle range calculated from time integration and DPSC-MC.}
\label{fig:ba_1d_distrib}
\end{figure}

To further investigate the convergence of the DPSC-MC approach, this section considers additional time integration cycles. To limit computational costs, this section uses the two-element model, which shows similar behavior with slight offsets in the absolute values compared to the 16-element case that was previously detailed. Using 10,000 stochastic time integration cycles instead of 1000 stochastic time integration cycles (both with 128 steps per cycle) results in the time integration covariance converging to the covariance calculated with DPSC-MC in \Cref{fig:ba_cov_dist_convg}. Furthermore, both approaches show a similar distribution that is clearly not Gaussian in \Cref{fig:ba_cov_dist_convg}. This non-Gaussian behavior is likely caused by evaluating the range of the bending moments, which is not linear for combinations of perturbations from the DPSC. Likewise, the marginal distributions in \Cref{fig:ba_1d_distrib} highlight the similarity in the sampled distributions of the fatigue response variables.
Alternatively, using 512 steps per cycle instead of 128 steps per cycle for time integration resulted in a similar improvement, even with only 1000 cycles of time integration. However, using smaller time steps can make the nonlinear solve more difficult because the perturbation size is inversely proportional to the length of the time step and thus grows for smaller time steps (see Appendix~\ref{sec:appendix_time_int}).

\FloatBarrier

\FloatBarrier

\FloatBarrier

\subsection{Numerical Conditioning} \label{sec:aero_conditioning}

\begin{figure}[h!]
\centering
        \includegraphics[scale=1,width=\linewidth]{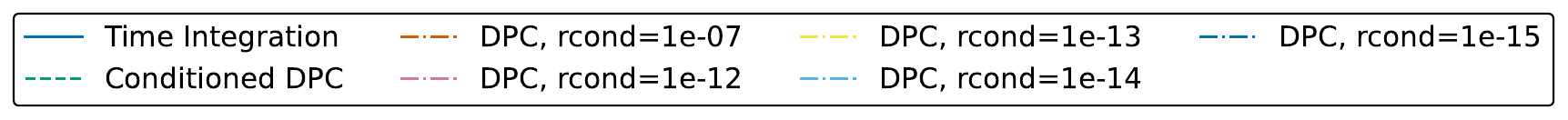} 
	\begin{subfigure}[b]{0.49\linewidth } 
		\centering
		\includegraphics[width=\linewidth]{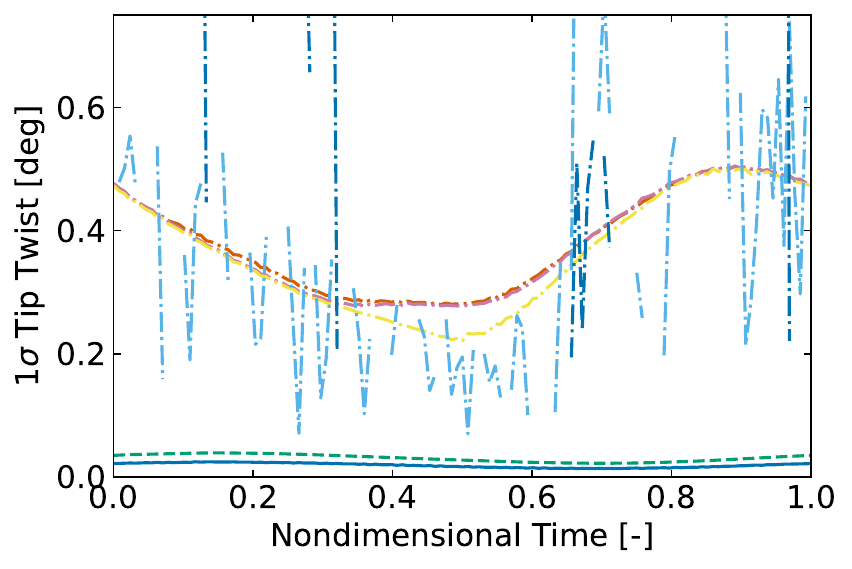}
		\caption{}
	\end{subfigure}	\hfill
	\begin{subfigure}[b]{0.49\linewidth } 
		\centering
		\includegraphics[width=\linewidth]{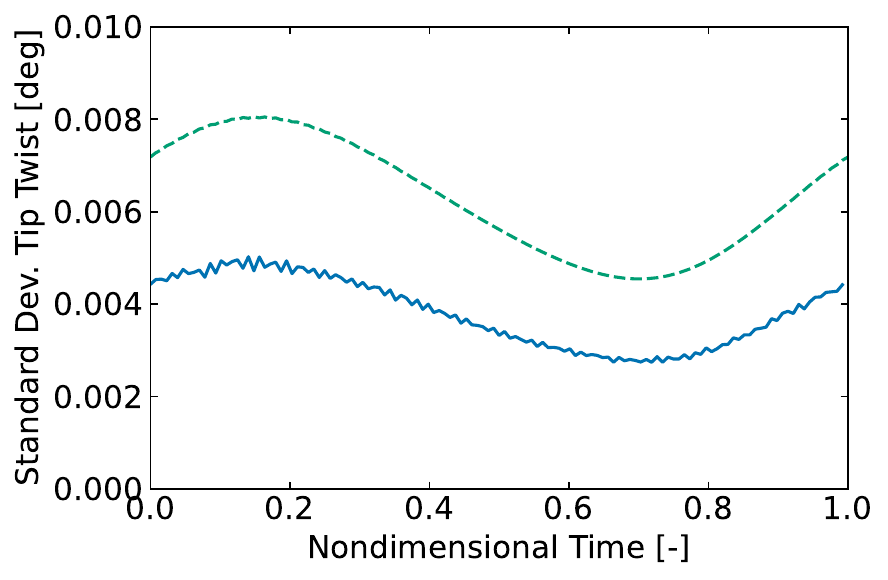}
		\caption{}
	\end{subfigure}\caption{Time varying standard deviation for tip twist motion of four-element model: (a) simulations with varying condition threshold for approximate solve and (b) zoomed in on time integration and reformulated solutions.}
\label{fig:ba_rcond}
\end{figure}

A key contribution of the present work is the formulation of the DPC method to avoid solving ill-conditioned linear systems.
A simpler approach is to use an approximate linear solve with a condition threshold (rcond in Python's NumPy least squares matrix solve\footnote{\url{https://numpy.org/doc/stable/reference/generated/numpy.linalg.lstsq.html}}). However, such an approach does not consistently work as shown in \Cref{fig:ba_rcond} for the four-element case.
Here, using rcond between 1e-7 and 1e-12 appears to give converged time-varying standard deviations for the tip twist response, but these differ significantly from the time integration results. Thus, that solution is still wrong, and it may be confused for a converged solution if only considering the DPC approach with the condition threshold. Further decreasing the conditioning threshold results in a slight decrease in the standard deviations before the solution goes numerically unstable and prevents the use of the results.
On the other hand, the conditioned DPC approach presented here does not suffer from these challenges.

\FloatBarrier

\FloatBarrier

In \Cref{fig:ba_rcond}, there is an offset between the conditioned DPC and time integration standard deviations.
Because the tip twist is much smaller than the responses in the other directions, this is not a major concern for this work. 
Time integration simulations with smaller time steps did show an increased standard deviation in the twisting at the tip for the two-element case compared to using the same 128 steps per cycle as used here. Thus, it is likely that the conditioned DPC approach represents a more accurate solution than the time integration solution in \Cref{fig:ba_rcond}. This is probably related to the fact that the first torsional natural frequency (about 4 Hz) is approximately 32 times the forcing frequency (0.127 Hz). Using 128 steps per cycle gives four steps per cycle of torsional mode perturbations and is likely insufficient to fully capture the behavior.

\section{Frictional System} \label{sec:friction_system}

This section applies the present methodology to a frictional system due to the importance of friction nonlinearities for many engineering systems \cite{brakeMechanicsJointedStructures2017}.
\Cref{sec:fric_describe} describes the system of a Euler-Bernoulli beam with a frictional tip contact based on \cite{balaji_quasi_static_2020}, and \Cref{sec:friction_deterministic} presents the deterministic response to harmonic forcing. 
The system is subjected to stochastic variations in the frictional slip force, and stochastic simulation results are shown in \Cref{sec:fric_stochastic_results} followed by a discussion of computational costs in \Cref{sec:fric_comp_costs}.
This system shows similar conditioning issues if the present modifications for the DPC are not applied as shown in \Cref{sec:fric_conditioning}.
Lastly, limitations in convergence are highlighted when the system becomes too nonlinear in \Cref{sec:fric_large_perturb}.

\subsection{System Description} \label{sec:fric_describe}

The frictional system uses seven cubic Euler-Bernoulli beam elements for bending in a single direction. The beam model parameters based on \cite{balaji_quasi_static_2020} are provided in \Cref{tab:friction_pars}. Rather than solving an axial load problem to determine the contact force, this work just uses a nominal 10 N frictional slip force.
Additionally, this work uses a Bouc-Wen model to more smoothly capture the stick-slip behavior.

\FloatBarrier

\begin{table}[h!]
\centering
\caption{Parameters for frictional system \cite{balaji_quasi_static_2020}.} \label{tab:friction_pars}
\begin{tabular}{clc}
\hline\hline
Symbol & Description    & Value                         \\ \hline
$l$ & Length & 1 m \\ 
$r$ & Radius (circular section) & 0.01 m \\
$E$ & Elastic modulus & 200 GPa \\
$\rho$ & Density & 7800 kg/m$^3$ \\
$\alpha_d$ & Mass proportional damping coefficient & 0.1 1/s \\
$\beta_d$ & Stiffness proportional damping coefficient & 1e-5 s \\
$f_s$ & Slip force & 10 N \\
$k_t$ & Tangential stiffness (at full slip reversal) & 10 kN/m \\
$\gamma$ & Ratio of Bouc-Wen parameters $\beta/\alpha$ & $-0.98$ \\
$n$ & Bouc-Wen exponent & 5 \\
\hline\hline
\end{tabular}
\end{table}

The Bouc-Wen model is a popular friction and hysteretic model \cite{fantetti_impact_2019, mathisReviewDampingModels2020, shetty_parametric_2023, porter_quantitative_2022}. Here, only a single frictional contact point is used and the frictional force is described as a first-order state $z$ with time derivative \cite{mathisReviewDampingModels2020}
\begin{equation}
    \dot{z} = \left[ A - \left[\alpha \text{sign}(z \dot{u})  - \beta\right] |z|^n \right] \dot{u},
\end{equation}
where $u$ is the displacement associated with the frictional element, $A$ is the stiffness at zero force, and the parameters $\alpha$, $\beta$, and $n$ control the shape of the hysteresis loop. Generally, $\alpha \geq 0$ and $\beta \in [-\alpha, \alpha]$ \cite{mathisReviewDampingModels2020}.
In this work, the Bouc-Wen model is used as a smoother alternative to a simple Jenkins element so that derivatives are well defined, but it is desired to have a relatively sharp transition between the stuck and slip regions.
To achieve this, a relatively high value of $n=5$ is used and the ratio $\gamma=\alpha/\beta=-0.98$ is chosen.
The values of $A$ and $\alpha$ can then be analytically calculated to match the desired slip force $f_s$ and tangential stiffness $k_t$ when reversing from full slipping. Specifically, this calculation is described by the equations
\begin{equation} \label{eq:bw_match_fskt}
\begin{split}
    A - (1 - \gamma) |f_s|^n \alpha &= 0 \\
    A + (1 + \gamma) |f_s|^n \alpha &= k_t.
\end{split}
\end{equation}
Here, the first condition is derived by requiring $\dot{z} = 0$ at the slip force, and the second condition is derived from having a tangential stiffness of $k_t$ when the velocity reverses while at the slipping force.
When $f_s$ is stochastically perturbed in this work, $A$ and $\alpha$ are updated so that \eqref{eq:bw_match_fskt} is still satisfied.

\subsection{Deterministic Response} \label{sec:friction_deterministic}

\begin{figure}[h!]
\centering
        \includegraphics[scale=1,width=0.5\linewidth]{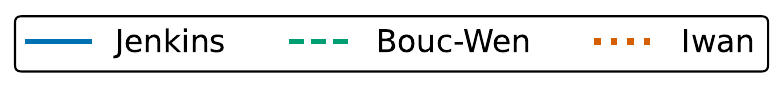} \\
	\begin{subfigure}[b]{0.49\linewidth } 
		\centering
        		\includegraphics[width=\linewidth]{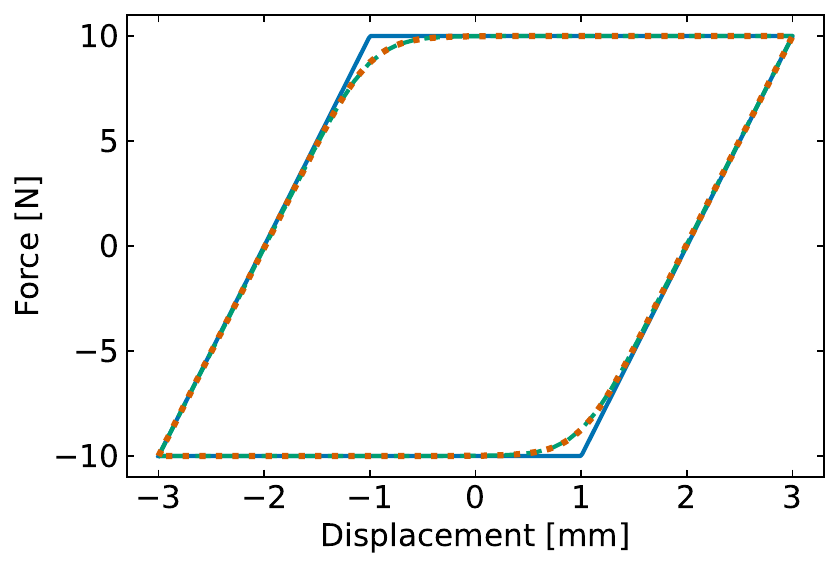}
		\caption{}
        \label{fig:fric_hyst_loops}
	\end{subfigure}	\hfill
	\begin{subfigure}[b]{0.49\linewidth } 
		\centering
        		\includegraphics[width=\linewidth]{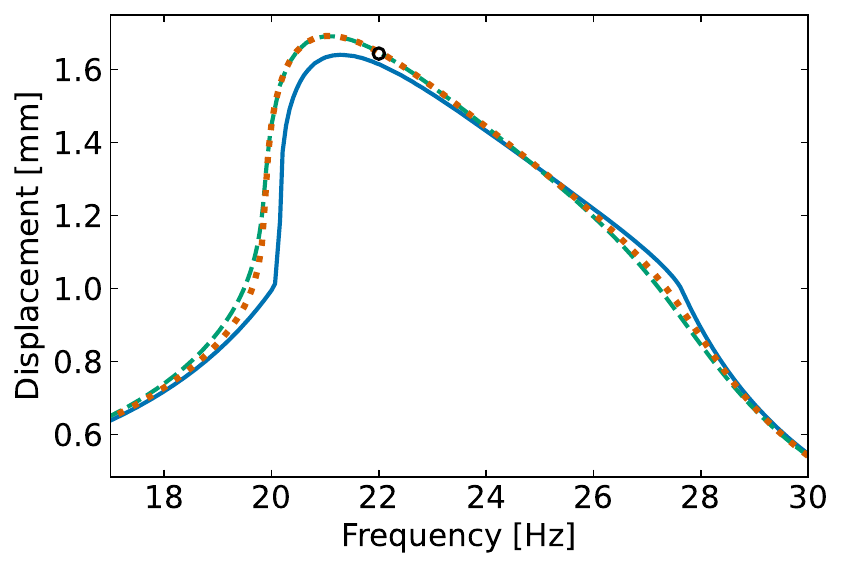}
		\caption{}
        \label{fig:fric_frc}
	\end{subfigure}\caption{Frictional system responses: (a) simply hysteresis loops and (b) frequency response curves. The circle in (b) denotes the single steady-state time integration solution for the Bouc-Wen model at the forcing frequency used in the subsequent discussion.}
\label{fig:friction_loop_frc}
\end{figure}

The hysteresis loop for the Bouc-Wen model is shown in comparison to a Jenkins element with the same slip force and tangential stiffness in \Cref{fig:fric_hyst_loops}. Using relatively high values of $n$ and $\gamma$ in the Bouc-Wen model allows it to closely match the Jenkins element except for in the transition region where the Bouc-Wen model shows the desired smoothing. 
An Iwan model is also derived from the response of the Bouc-Wen model going from full negative to positive slip\footnote{A total displacement range of $4 f_s / k_t$ is used with the force starting at $-f_s$.} and matches the Bouc-Wen hysteresis loop in \Cref{fig:fric_hyst_loops}.
The Iwan model is a sum of Jenkins elements with stiffness and slip forces based on the Bouc-Wen hysteresis loop rather than a standard distribution of sliders.

\FloatBarrier

\FloatBarrier

\Cref{fig:fric_frc} shows the frequency response curves (FRCs) for the three frictional models subject to a harmonic tip force with magnitude 5~N. 
Here, the HBM solutions for the FRCs use the zeroth and first three harmonics. The AFT procedure and subsequent DPC calculations use 128 steps per period.
The gradients from the Iwan model are used in the HBM solve for the FRC of the Bouc-Wen model.
The gradients of the Iwan model are significantly easier to implement in the AFT scheme, so the present work exploits the available implementation of those gradients in TMDSimPy \cite{porterTMDSimPy}.
However, the Bouc-Wen FRC uses the force from the Bouc-Wen model in the residual calculation, so it accurately captures the behavior of the system with the Bouc-Wen frictional element. 
The use of the gradients from the Iwan model only affects the convergence of the model, but given the similarity in the hysteresis loops, no convergence issues were noted.
Since the Bouc-Wen model is not a sum of simple stick-slip elements, up to 30 cycles of the hysteresis loop are calculated using a trapezoid time integration scheme to obtain a steady-state response used in the HBM residual calculation.
This was sufficient to close the hysteresis loop of the Bouc-Wen model to within a tolerance of 1e-4 relative to the maximum force anywhere in the cycle for the entire FRC shown in \Cref{fig:fric_frc}. Between approximately 18.5 Hz and 28 Hz, the hysteresis loops of the Bouc-Wen model were closed to within a tolerance of 1e-10.

The FRCs are typical for a frictional system (e.g., see \cite{balaji_quasi_static_2020}). The Bouc-Wen and Iwan FRCs deviate from the Jenkins element, but may still be sufficiently close if one desires to change to the smoother model to use the subsequent stochastic analysis. Furthermore, smooth stick-slip transitions are generally considered more realistic for frictional systems \cite{porter_towards_2023, pesaresi_advanced_2018, porter_quantitative_2022}. To more closely match the Bouc-Wen model to the Jenkins element, either a larger value for $n$ or a smaller value for $\gamma \in [-1, 1]$ would be required. However, further increasing $n$ or decreasing $\gamma$ can introduce numerical stability or conditioning issues.
As expected, the Iwan model that was fit to the Bouc-Wen hysteresis loop behaves similarly to the Bouc-Wen model with slight deviations at around 1 mm amplitudes in the transition to macroslip regime. This is because the Bouc-Wen model does not match the tangential stiffness $k_t$ when reversing from forces less than $f_s$, whereas the Iwan model always has tangential stiffness $k_t$. Despite this difference in the frictional force, the FRCs remain similar.

The time integration solution at 22 Hz for the Bouc-Wen model is consistent with the HBM solution in \Cref{fig:fric_frc}. 
The time integration solution uses 128 time steps per excitation cycle and reached steady state after 15 cycles of time integration.
That point is further investigated in the next section when incorporating stochastic excitations.

\subsection{Stochastic Response} \label{sec:fric_stochastic_results}

\begin{figure}[h!]
    \centering
            
	\begin{subfigure}[b]{0.49\linewidth } 
		\centering
                		\includegraphics[width=\linewidth]{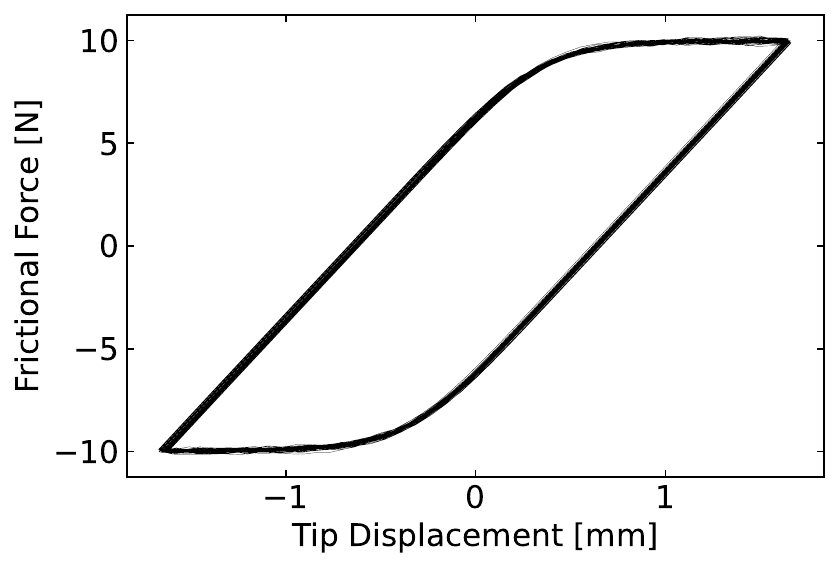}
		\caption{}
	\end{subfigure}	\hfill
	\begin{subfigure}[b]{0.49\linewidth } 
		\centering
                		\includegraphics[width=\linewidth]{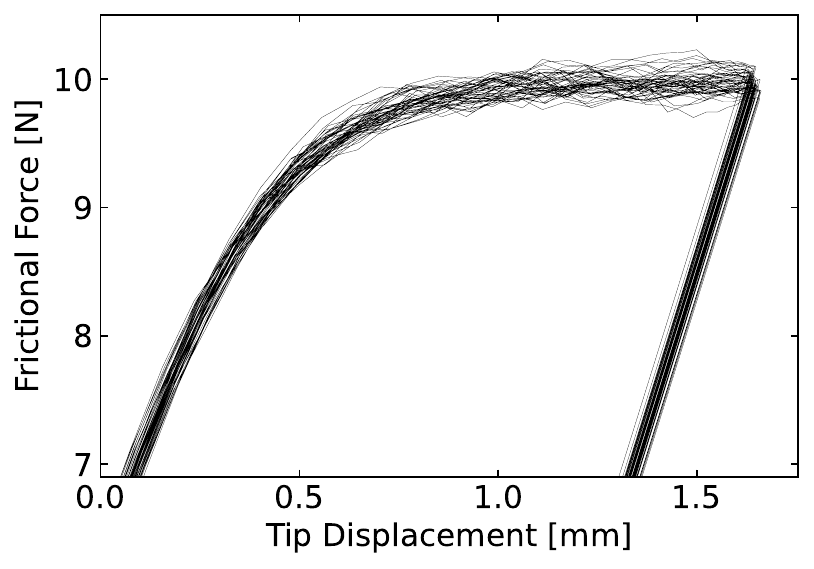}
		\caption{}
	\end{subfigure}    \caption{The last 50 hysteresis loops calculated via time integration for the stochastically excited frictional system. Plot (b) is a zoomed-in portion of (a).}
    \label{fig:bw_stochastic_hyst}
\end{figure}

For the stochastic analysis, the frictional system is subject to a harmonic tip load of 5~N at 22~Hz and the slip force $f_s$ varies with a standard deviation of $\sigma = 0.004$~N. The noise in the slip force could represent variations in the normal force that holds two surfaces together or in the contact surfaces that are slipping against each other.
\Cref{fig:bw_stochastic_hyst} provides further context on the effect of the stochastic variations in the slip force. While slipping, the variances in the slip force cause deviations in the hysteresis loops. In the stuck regions after reversal, each hysteresis loop appears linear but offset due to the variations immediately before reversal.

\FloatBarrier

\begin{figure}[h!]
    \centering
        \includegraphics[width=0.7\linewidth]{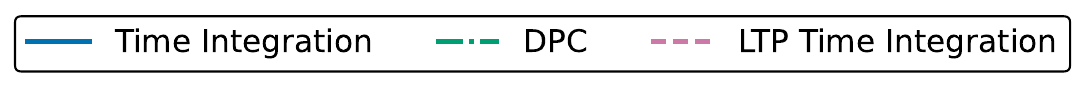}
	\begin{subfigure}[b]{0.49\linewidth } 
		\centering
                		\includegraphics[width=\linewidth]{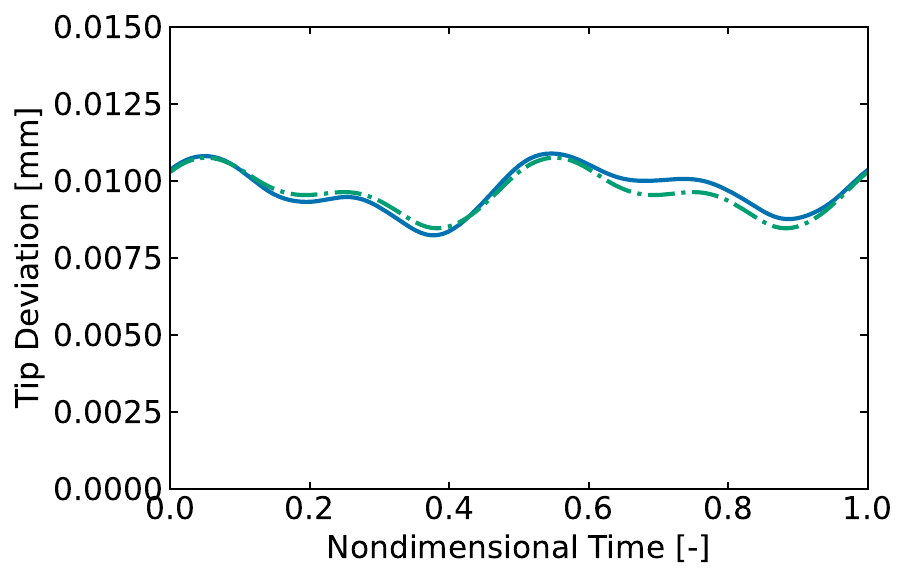}
		\caption{}
        \label{fig:bw_std}
	\end{subfigure}	\hfill
	\begin{subfigure}[b]{0.49\linewidth } 
		\centering
                		\includegraphics[width=\linewidth]{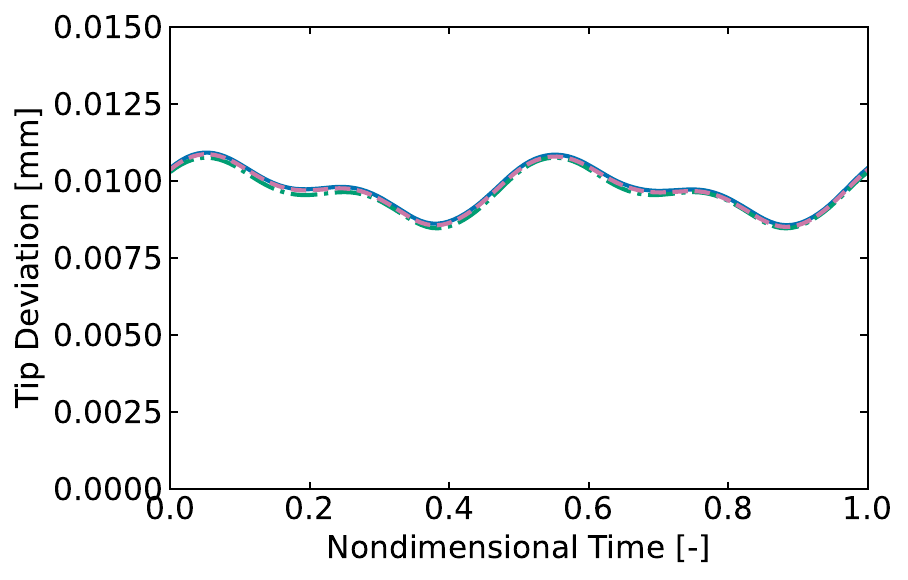}
		\caption{}
        \label{fig:bw_std_ltp_convg}
	\end{subfigure}\caption{Time-varying standard deviation at the tip using (a) 1000 cycles of time integration and (b) 10,000 cycles of time integration.}
\label{fig:bw_stochastic_simple}
\end{figure}

The time-varying standard deviation of the tip displacement is shown in \Cref{fig:bw_stochastic_simple}. Clearly, the white noise variations in the slip force result in different amounts of variation for different points in the cycle. 
\Cref{fig:bw_std} shows some differences between the time-varying standard deviations from time integration and DPC. These differences are attributed to the 1000 cycles at 128 steps per cycle for time integration being insufficient to fully converge the response statistics.
This is confirmed in \Cref{fig:bw_std_ltp_convg} where 10,000 cycles of time integration give a standard deviation consistent with DPC. 
Alternatively, increasing to 256 steps per cycle for 1000 cycles (not shown) improved the match between the time integration DPC solutions compared to 128 steps per cycle for 1000 cycles.
Also in \Cref{fig:bw_std_ltp_convg}, time integration is done on the LTP system used in DPC and based on the linearization around the HBM solution. This LTP system is useful for isolating convergence differences from the unaccounted-for effects of nonlinearity in DPC in subsequent analysis.

\begin{figure}[h!]
    \centering
        \includegraphics[width=\linewidth]{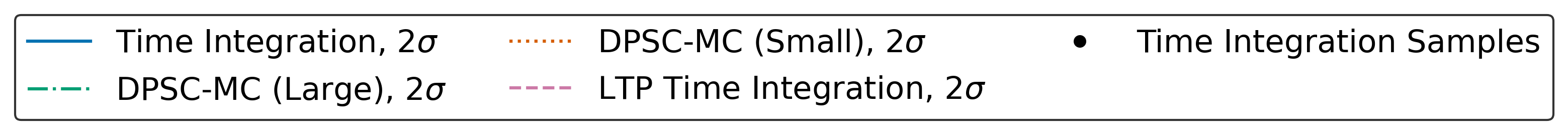}
    	\begin{subfigure}[b]{0.49\linewidth } 
		\centering
                		\includegraphics[width=\linewidth]{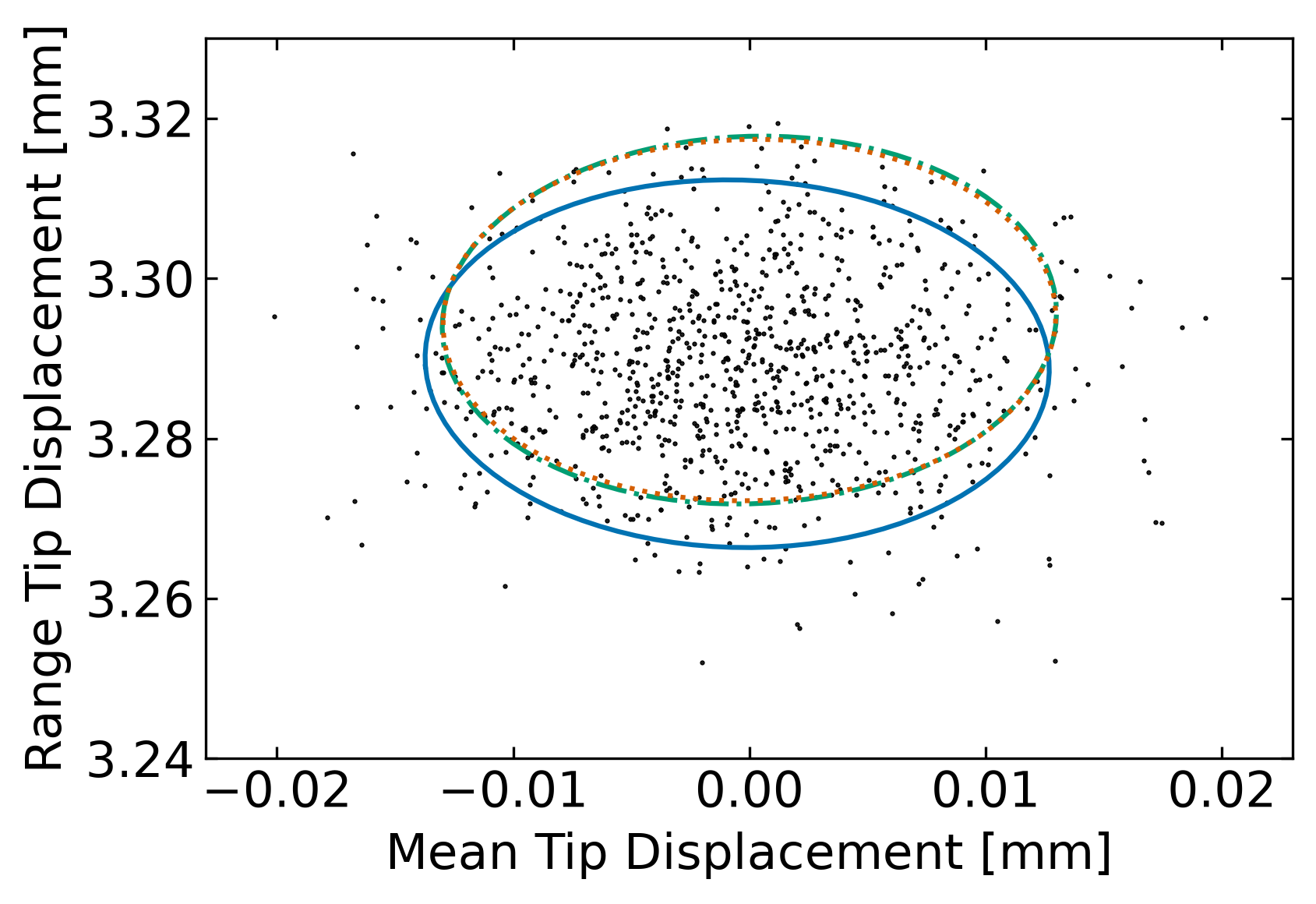}
		\caption{}
	\end{subfigure}	\hfill
	\begin{subfigure}[b]{0.49\linewidth } 
		\centering
                		\includegraphics[width=\linewidth]{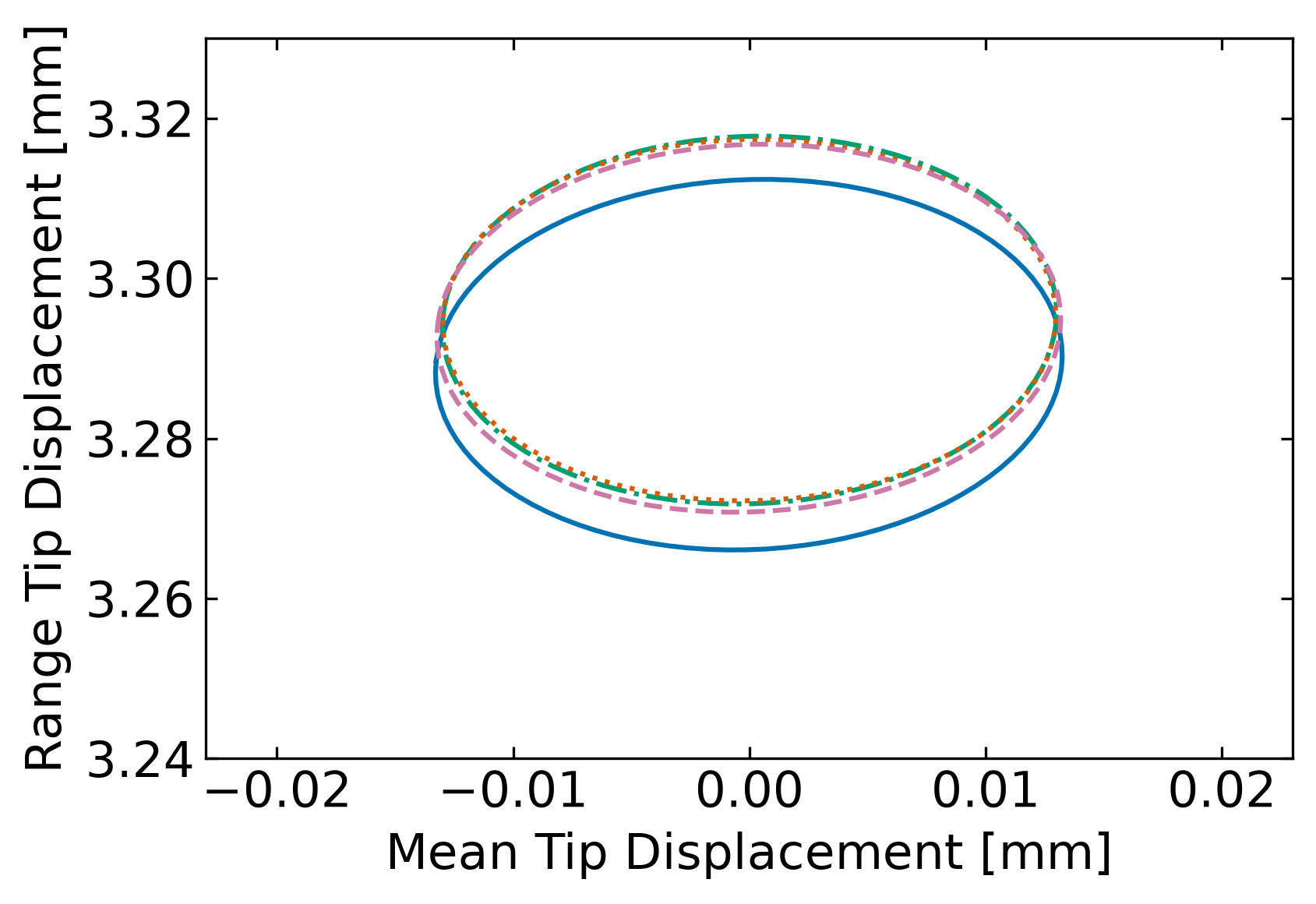}
		\caption{}
        \label{fig:bw_cov_convg}
	\end{subfigure}    \caption{Comparison of covariances calculated via time integration and DPSC-MC: (a) 1000 cycles of time integration are used to calculate the covariance and are plotted as individual points; (b) 10,000 cycles of time integration are used for both the baseline time integration and the LTP time integration solutions.}
    \label{fig:bw_covariance}
\end{figure}

Building on the DPC, the DPSC-MC results are shown in \Cref{fig:bw_covariance}. The DPSC-MC approach gives a clear offset compared to the time integration solution for the mean value across stochastic cycles of the range of tip displacements within a cycle. \Cref{fig:bw_cov_convg} highlights that this difference is not related to convergence. Rather, DPSC-MC matches the LTP time integration covariance, but the errors in both compared to the baseline time integration are attributed to the lack of nonlinear effects in DPSC-MC and the LTP time integration solutions.
However, the DPSC-MC solution is significantly faster than the base time integration solution (see \Cref{sec:fric_comp_costs}), so the results still provide a faster useful estimate despite some errors in the covariance calculation.
It is also noted that these errors in the DPSC-MC calculation occur despite no clear errors in the DPC in \Cref{fig:bw_stochastic_simple}. This implies that the present nonlinear effects are specifically important for how different time instants are correlated.

\FloatBarrier

\FloatBarrier

\FloatBarrier

\subsection{Computation Time} \label{sec:fric_comp_costs}
\newcommand{\bwSteadyNumCycles}[0]{15}
\newcommand{\bwSteadyTimeSec}[0]{3.9}
\newcommand{\bwStochasticTimeSec}[0]{711.5}
\newcommand{\bwStochasticCycles}[0]{1000}
\newcommand{\bwContinuationTime}[0]{36.1}
\newcommand{\bwContinuationNumSteps}[0]{14}
\newcommand{\bwSingleHBMTime}[0]{2.0}
\newcommand{\bwFullDirectTime}[0]{28.6}
\newcommand{\bwMinDirectTime}[0]{3.2}
\newcommand{\bwFullDirectSpeedup}[0]{24.8}
\newcommand{\bwMinDirectSpeedup}[0]{222}
\newcommand{\bwStochasticLTPTimeSec}[0]{1101.0}
\newcommand{\bwStochasticLTPCycles}[0]{10,000}

As with the aeroelastic examples, there is a significant computational benefit to applying the DPSC-MC compared to stochastic time integration.
Stochastic time integration with {\bwStochasticCycles} cycles took {\bwStochasticTimeSec} s.
On the other hand, the direct method including convergence checking for the DPSC and the large Monte Carlo sampling took {\bwFullDirectTime} s, representing a {\bwFullDirectSpeedup}$\times$ speedup.
However, this speedup would be larger if compared to also checking convergence for the statistics calculated from time integration.
During design, convergence checks may not be necessary at every iteration, so without the convergence checks and using the small Monte Carlo set, DPSC-MC took {\bwMinDirectTime} s for a computational speedup of {\bwMinDirectSpeedup}$\times$. This is approximately equivalent to obtaining the stochastic statistics in the computation time equivalent to the time integration of 4.5 stochastic cycles.

For the frictional system, HBM failed to converge when attempting to solve the solution directly at the point of interest from a linear prediction with no stiffness from friction. 
Rather, this work used a continuation from 23 Hz to 21 Hz to identify HBM solutions taking {\bwContinuationNumSteps} steps in {\bwContinuationTime} s. 
Using the solution from the frequency closest to the desired 22 Hz as the initial guess, a single HBM solve took {\bwSingleHBMTime} s. Overall, this is not a clear speedup compared to the computation time of {\bwSteadyTimeSec} s for time integration to calculate {\bwSteadyNumCycles} cycles and converge to steady state.
The lack of a significant speedup for HBM can be partially attributed to the fact that the AFT algorithm relies on time integration of multiple cycles of the Bouc-Wen model to reach steady state and the relatively small problem size evaluated here.

The computation time for DPSC-MC benefits from only considering a linearized system rather than the full nonlinear system. Alternatively, similar statistics are obtained using the LTP time integration formulation.
The LTP time integration solution is faster than the full nonlinear stochastic time integration and took {\bwStochasticLTPTimeSec} s for {\bwStochasticLTPCycles} cycles. Scaling this to a reference of 1000 cycles of time integration with the LTP system results in a computation time of 110.1 s. Compared to this estimated LTP computation time, the small DPSC-MC has a speedup of 34.4$\times$. The LTP timing here does not include the creation of the state-space model that is used for both the DPSC-MC and LTP time integration calculations and represents about half of the computation time for the small DPSC-MC calculation.

\subsection{Numerical Conditioning} \label{sec:fric_conditioning}

\begin{figure}[h!]
\centering
        \includegraphics[width=0.75\linewidth]{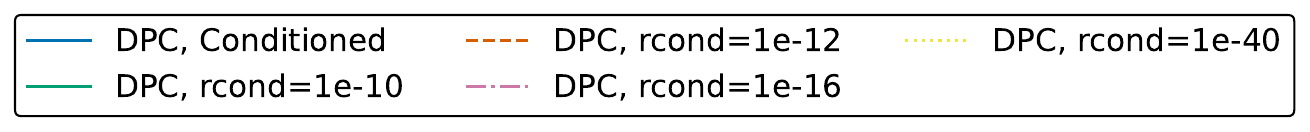}\\
	\begin{subfigure}[b]{0.49\linewidth } 
		\centering
        		\includegraphics[width=\linewidth]{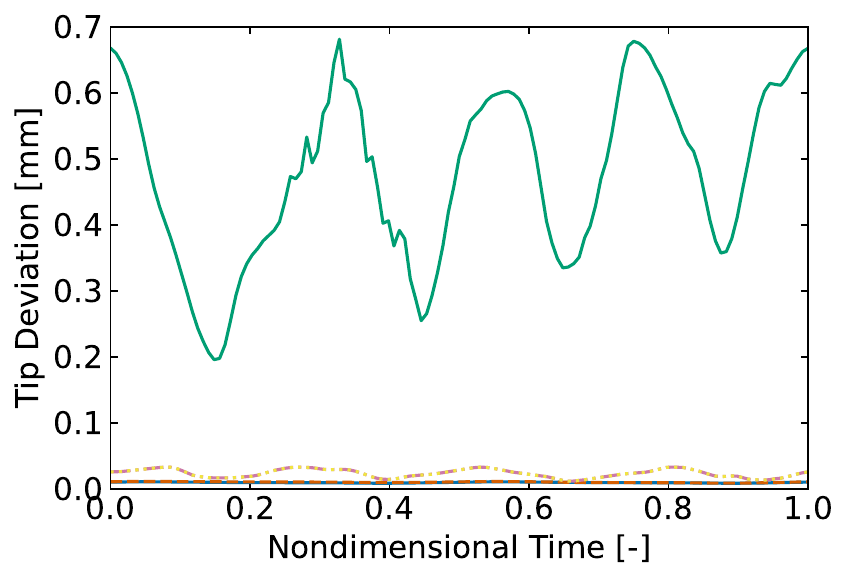}
	\end{subfigure}	\hfill
	\begin{subfigure}[b]{0.49\linewidth } 
		\centering
        		\includegraphics[width=\linewidth]{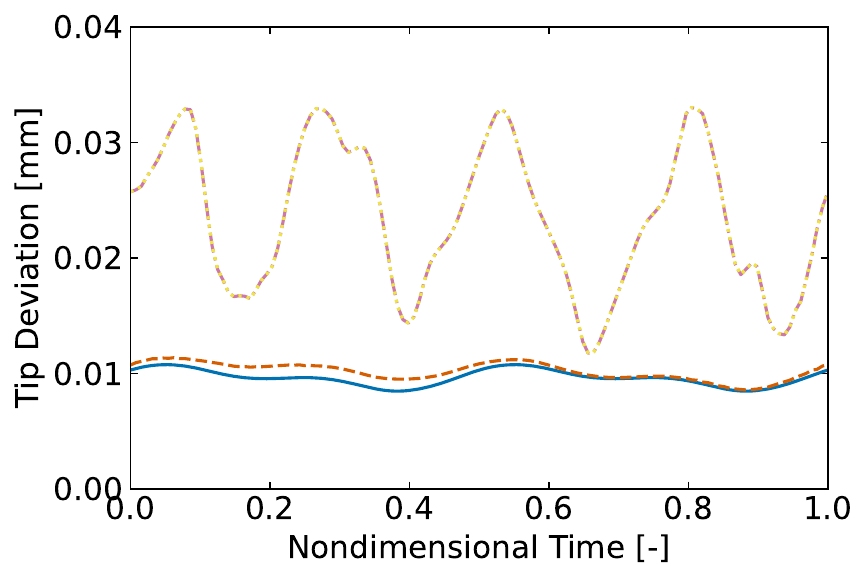}
	\end{subfigure}\caption{Time-varying standard deviation with different conditioning thresholds for the 21 element model. Plots are the same data with different $y$-axis scales.}
\label{fig:bw_rcond}
\end{figure}

Given the conditioning challenges for the aeroelastic system, potential conditioning challenges were investigated for the frictional system.
Using 7 elements, the conditioning threshold rcond could be adjusted to converge the time-varying standard deviation to that obtained through the conditioned DPC approach.\footnote{As shown previously, the time integration solutions converge to the conditioned DPC solution.}
However, \Cref{fig:bw_rcond} shows that for 21 elements, adjusting rcond is insufficient to match the conditioned solution.
While the rcond threshold of 1e-12 was roughly close to the conditioned DPC solution, further decreasing the threshold resulted in the solution diverging.
This highlights the importance of using the conditioned algorithm. Furthermore, since the issues appear when increasing from 7 to 21 elements, the conditioned algorithm could be more important for even larger models at engineering scale.

\FloatBarrier

\FloatBarrier

\subsection{Larger Deviations} \label{sec:fric_large_perturb}

\begin{figure}[h!]
    \centering
        
    \begin{subfigure}[b]{0.49\linewidth } 
		\centering
                		\includegraphics[width=\linewidth]{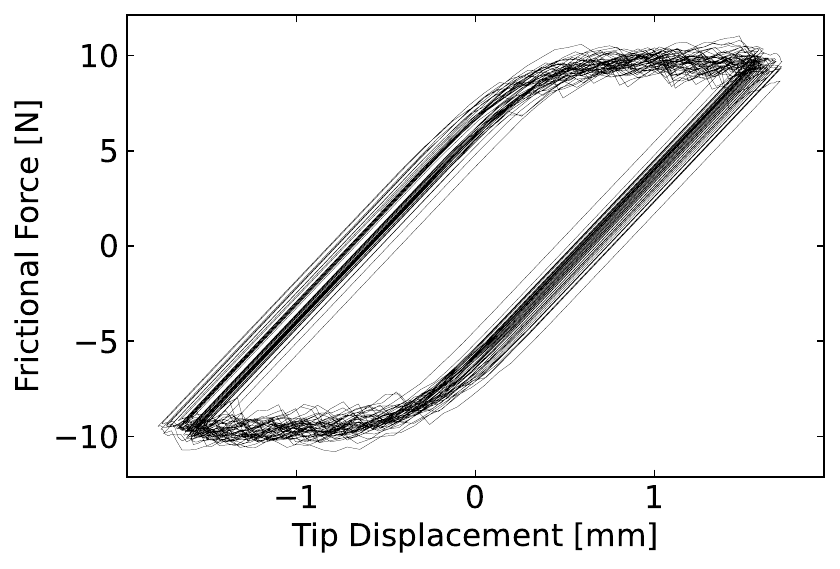}
		\caption{}
	\end{subfigure}	\hfill
	\begin{subfigure}[b]{0.49\linewidth } 
		\centering
                		\includegraphics[width=\linewidth]{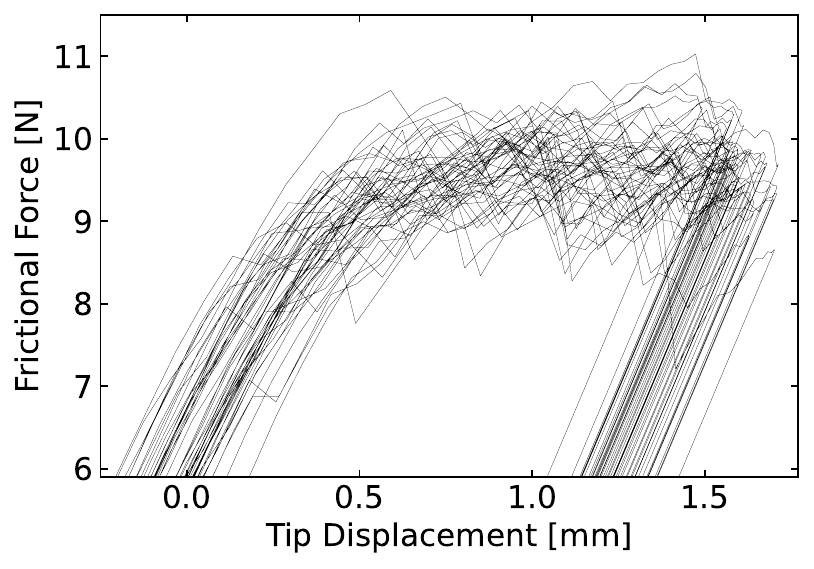}
		\caption{}
	\end{subfigure}    \caption{Time-varying responses for the last 50 hysteresis loops calculated via time integration for the stochastic response of the Bouc-Wen and beam model. Plot (b) is a zoomed-in portion of (a).}
    \label{fig:bw_stochastic_hyst_largeSigma}
\end{figure}

\begin{figure}[h!]
    \centering
            \includegraphics[width=0.7\linewidth]{bw_tip_std_case1_LTP_legend.pdf}
	\begin{subfigure}[b]{0.49\linewidth } 
		\centering
                		\includegraphics[width=\linewidth]{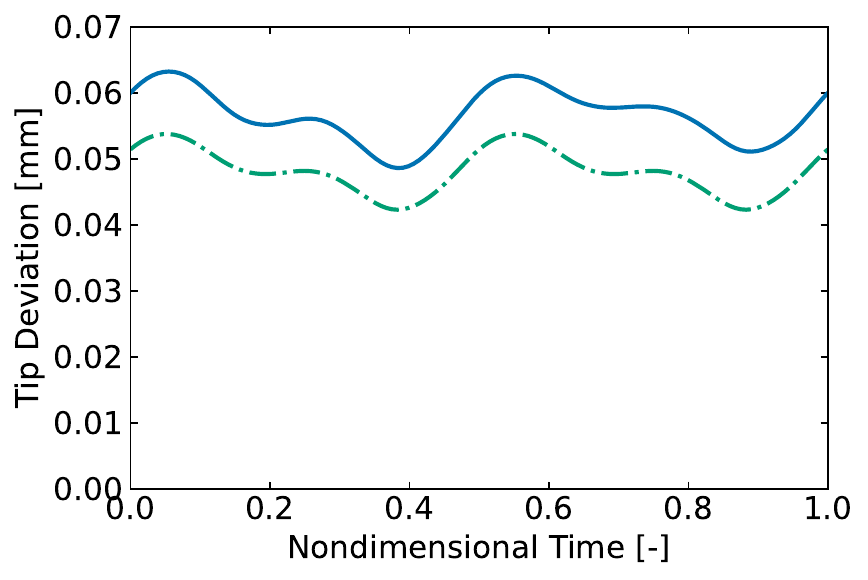}
		\caption{}
        \label{fig:bw_std_baseline_largeSigma}
	\end{subfigure}	\hfill
	\begin{subfigure}[b]{0.49\linewidth } 
		\centering
                		\includegraphics[width=\linewidth]{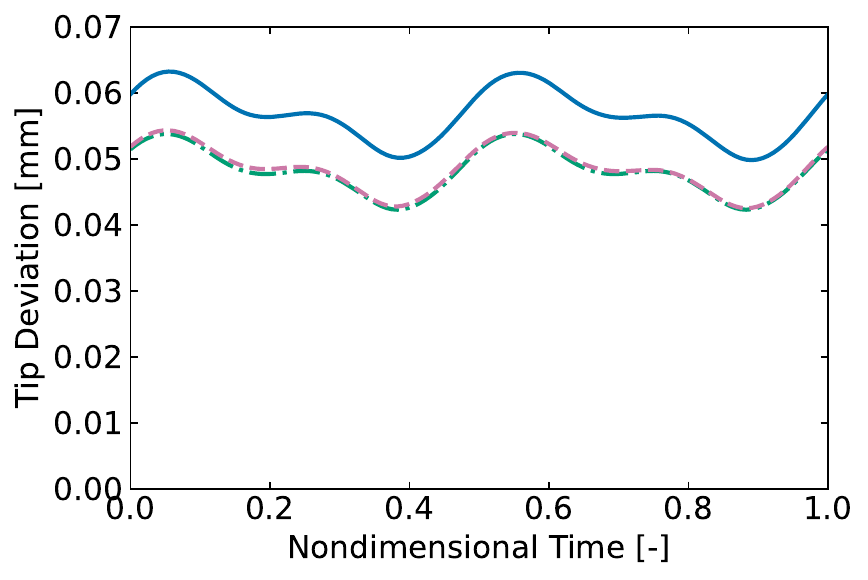}
		\caption{}
        \label{fig:bw_std_ltp_convg_largeSigma}
	\end{subfigure}\caption{Time-varying standard deviation at the tip: (a) 1000 cycles of time integration are used; (b) 10,000 cycles of time integration are used for both the baseline and LTP cases.}
\label{fig:bw_std_largeSigma}
\end{figure}

\begin{figure}[h!]
    \centering
        \includegraphics[width=\linewidth]{bw_tip_cov_case1_LTP_legend.png}
    	\begin{subfigure}[b]{0.49\linewidth } 
		\centering
                		\includegraphics[width=\linewidth]{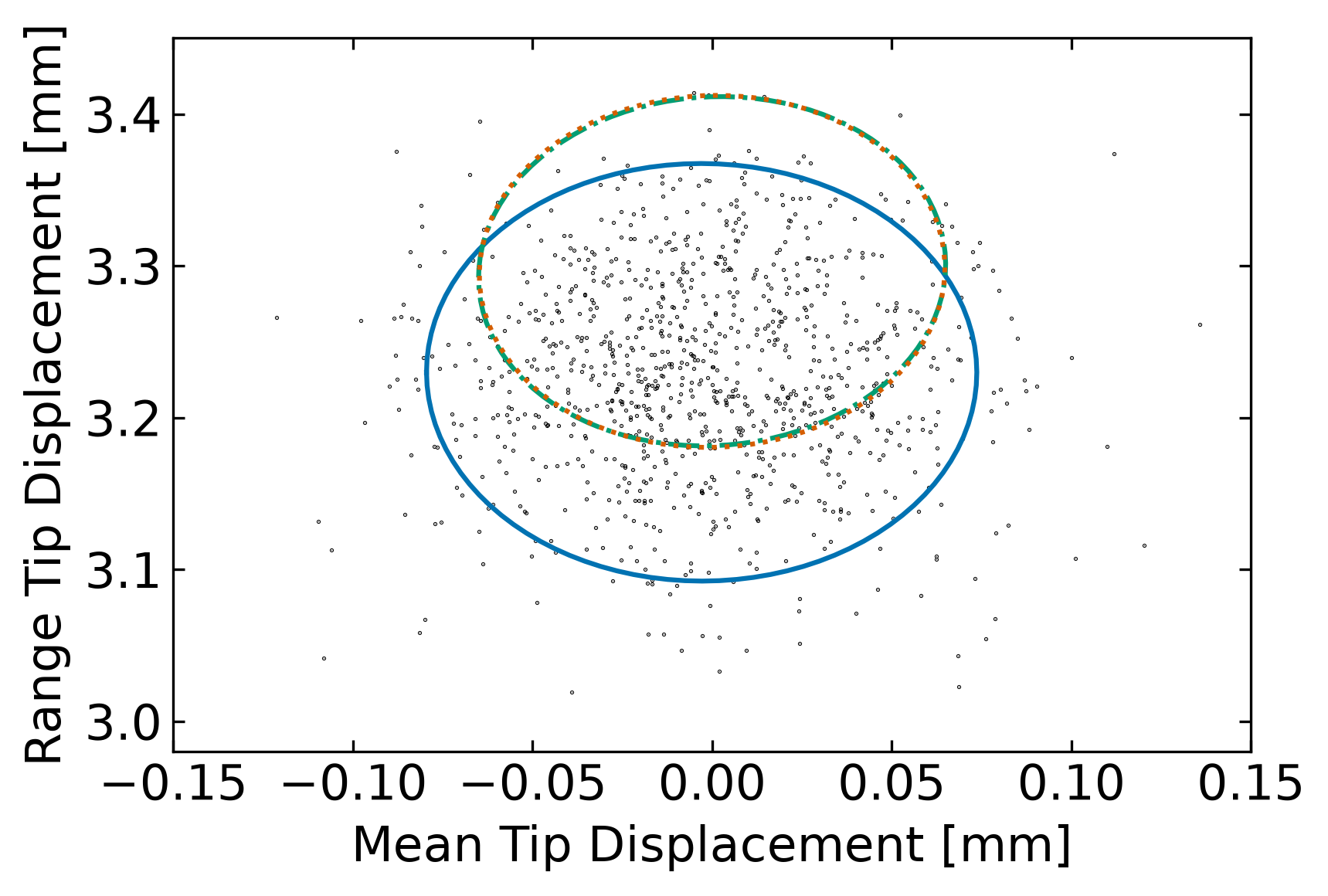}
		\caption{}
        \label{fig:bw_cov_1_largeSigma}
	\end{subfigure}	\hfill
	\begin{subfigure}[b]{0.49\linewidth } 
		\centering
                		\includegraphics[width=\linewidth]{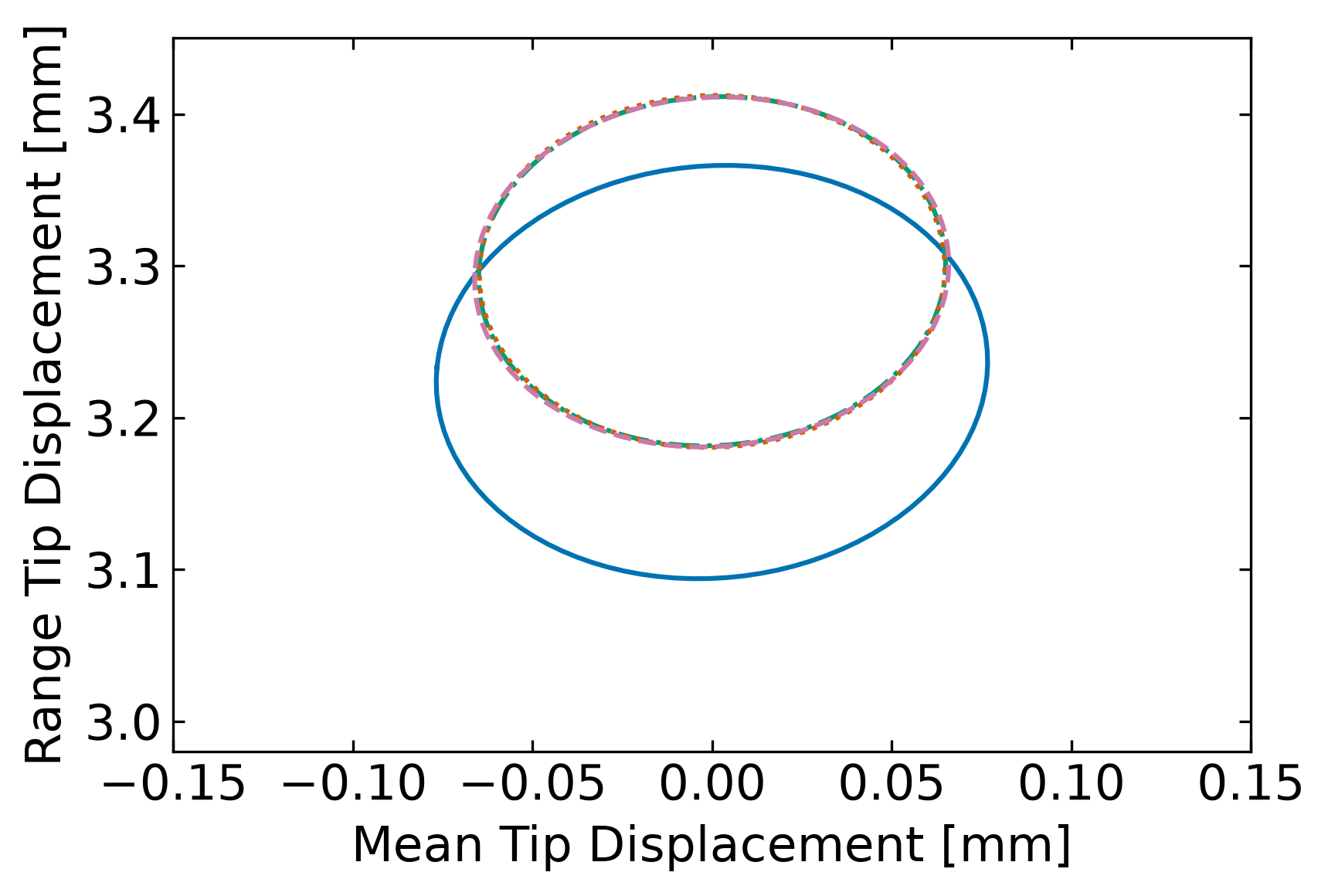}
		\caption{}
        \label{fig:bw_cov_2_largeSigma}
	\end{subfigure}    \caption{Comparison of covariances calculated via time integration and direct method. In (a), 1000 cycles of time integration are used to calculate the covariance and are plotted as individual points. In (b), 10,000 cycles of time integration are used for both the baseline time integration and the LTP time integration solutions.}
    \label{fig:bw_covariance_large_sigma}
\end{figure}

This section considers the limitations of the DPC and DPSC-MC approaches for larger deviations of $\sigma = 0.02$~N for the slip force.
This standard deviation significantly increases the deviation between the hysteresis loops in \Cref{fig:bw_stochastic_hyst_largeSigma} compared to \Cref{sec:fric_stochastic_results}.
In \Cref{fig:bw_std_largeSigma}, the larger deviation results in the DPC approach notably underestimating the time integration time-varying standard deviations. 
As before, this is attributed to the linearization simplification in DPC since DPC is consistent with LTP time integration in \Cref{fig:bw_std_ltp_convg_largeSigma}.
The errors in the time-varying standard deviation from DPC also propagate to the covariance of fatigue statistics in \Cref{fig:bw_covariance_large_sigma}. Again, these errors are attributed to the linearization simplifications in the present approach since DPSC-MC gives consistent results with stochastic time integration of the LTP system.
Despite these errors, the significant computational savings for DPSC-MC (see \Cref{sec:fric_comp_costs}) could mean that it is useful for approximating fatigue statistics.

\FloatBarrier

\FloatBarrier

\section{Conclusions} \label{sec:conclusions}

This work analyzes the effects of combined stochastic and periodic excitations. Responses to periodic excitations are deterministically calculated, and then linearized perturbations for the periodic system are analyzed by extending the methodology proposed in \cite{ahsan_adjoint_based_2025}.
The method of \cite{ahsan_adjoint_based_2025} provides for a direct calculation of a time-varying covariance of the response, and this work focuses on scaling the methodology to more complicated and larger systems.
The following are the specific contributions of this work:
\begin{itemize}
    \item The algorithm from \cite{ahsan_adjoint_based_2025} for the direct periodic covariance is modified to improve numerical conditioning. The importance of this modification is highlighted through two examples.
    \item The algorithm is extended to analyze the covariance of responses that span both spatial dimensions and different time instants in the period with a period-spanning covariance that uses a proposed reduction step to maintain computational efficiency.
    \item The period-spanning covariance is used to directly sample full stochastic cycles and rapidly estimate fatigue-relevant statistics.
    \item An aeroelastic system with white noise variations in the inflow air speed and a frictional system with white noise variations in the slipping force are analyzed.
    \item For the aeroelastic system, fatigue statistics are accurately estimated at computational speedups of {\stochShortSpeedupMin} to {\stochShortSpeedupMax} compared to stochastic time integration.
    \item For the frictional system, statistics are estimated at a speed-up of {\bwMinDirectSpeedup}$\times$, but the linearization of the present approaches introduces some errors compared to fatigue statistics from stochastic time integration.
\end{itemize}
Future work could analyze the effect of other types of noise by including additional states to filter the white noise \cite{spanos_filter_1986, alotta_fractional_2014, breunung_noise_2023}.
Additionally, there is an opportunity to further extend the present work to also consider self-excited or quasi-periodic systems that are also considered in \cite{ahsan_adjoint_based_2025}.
For the self-excited case from \cite{ahsan_adjoint_based_2025}, it is likely that the present conditioned modifications would also be necessary. A possible extension of the present conditioned algorithm to the self-excited case is briefly outlined in Appendix~\ref{sec:app_conditioned_self_excited} but is not investigated here.

Overall, this work provides for fast estimation of response statistics due to stochastic excitation near a nonlinear deterministic response. This provides a basis for expanding common nonlinear analyses beyond deterministic responses without relying on costly time integration.
The present approach is straightforward to apply for many nonlinear simulations because the methodology relies on gradient information that is already generally available as it is used in stability analyses.
Code used in the present work and numerical solution data are made available at~\cite{porterStochasticCodeData2026}.

\appendix

\section{Solution Methods} \label{sec:appendix_methods}

\Cref{sec:appendix_time_int} provides an overview of the deterministic and stochastic time integration methods used in this work.
Then, \Cref{sec:appendix_methods_hbm} reviews the harmonic balance method.

\subsection{Time Integration} \label{sec:appendix_time_int}

Time integration is conducted based on the generalized-$\alpha$ algorithm \cite{hughesFiniteElementMethod2012} for the second-order states $\vect{y}$  and with the trapezoid method for the first-order states $\vect{z}$. 
At each time integration step, all states are known at the starting time $t^n$ and are to be calculated at time $t^{n+1}$ for the time step of length $\Delta t = t^{n+1} - t^n$.
The equations of motion for the first- and second-order states are given by \eqref{eq:1stand2ndODE}.

The algorithm for the second-order states uses parameter $\alpha \in [-^1/_3, 0]$ (chosen here as 0) and derived parameters
\begin{equation}
\begin{split}
    \beta &= \dfrac{(1 - \alpha)^2}{4} \\
    \eta &= \dfrac{1 - 2 \alpha}{2} .
\end{split}
\end{equation}
During a time step, the accelerations $\vect{\ddot{y}}^{n+1}$ and first-order states $\vect{z}^{n+1}$ at the end of the time step are assumed and iteratively solved for.
Displacement and velocities at the end of the time step are calculated as
\begin{equation}
\begin{split}
    \vect{y}^{n+1} &= \vect{y}^n 
        + \Delta t \vect{\dot{y}}^n 
        + \dfrac{\Delta t^2}{2} \left[ (1 - 2 \beta) \vect{\ddot{y}}^n + 2 \beta \vect{\ddot{y}}^{n+1} \right]
    \\
    \vect{\dot{y}}^{n+1} &= \vect{\dot{y}}^n + \Delta t \left[ (1 - \eta) \vect{\ddot{y}}^n + \eta \vect{\ddot{y}}^{n+1} \right].
\end{split}
\end{equation}
The residual associated with the second-order states is
\begin{equation} \label{eq:2ndOrderTimeRes}
\begin{split}
    \vect{R_{2nd}} 
    &= \mat{M}\ \vect{\ddot{y}}^{n+1}
    + \mat{C}\left[ (1 + \alpha) \vect{\dot{y}}^{n+1} - \alpha \vect{\dot{y}}^n \right]
    + \mat{K}\left[ (1 + \alpha) \vect{y}^{n+1} - \alpha \vect{y}^n \right]
    \\ & \ \ \ +
    (1 + \alpha) \vect{f_{nl}}(\vect{y}^{n+1}, \vect{\dot{y}}^{n+1}, \vect{z}^{n+1}, \vect{p}, t^{n+1}) 
    - \alpha \vect{f_{nl}}(\vect{y}^n, \vect{\dot{y}}^n, \vect{z}^n, \vect{p}, t^n)
    \\ & \ \ \ -
    (1 + \alpha) \vect{f_{ext}}(t^{n+1}) 
    + \alpha \vect{f_{ext}}(t^n),
\end{split}
\end{equation}
where $\vect{f_{ext}}(t)$ is the external deterministic force.
The residual for the first-order states is
\begin{equation} \label{eq:1stOrderTimeRes}
    \vect{R_{1st}} = \vect{z}^{n} 
    + \dfrac{\Delta t}{2} \left[ 
        \vect{g}(\vect{y}^n, \vect{\dot{y}}^n, \vect{z}^n, \vect{p}, t^n)
        +
        \vect{g}(\vect{y}^{n+1}, \vect{\dot{y}}^{n+1}, \vect{z}^{n+1}, \vect{p}, t^{n+1})
    \right]
    - 
    \vect{z}^{n+1}.
\end{equation}
At each time step, \eqref{eq:2ndOrderTimeRes} and \eqref{eq:1stOrderTimeRes} are solved simultaneously with a monolithic Newton-Raphson algorithm \cite{porterTMDSimPy}.

A standard approach for stochastic time integration for first-order differential equations is the Euler-Maruyama scheme with time integration updates of \cite{kloeden_numerical_1999}
\begin{equation} \label{eq:apply_time_perturb}
    \vect{x}^{n+1} 
    = 
    \vect{x}^n 
    + \vect{\dot{x}}^n \Delta t 
    + \left[ 
        \dfrac{\sigma}{\sqrt{\Delta t}} 
        \mat{F}\left(\vect{x}^n\right) 
        \vect{\Delta W}^{n+1/2}
    \right]\Delta t.
\end{equation} 
Here, $\sigma$ is the standard deviation of the noise, and the noise is distributed according to the $ N \times m$ matrix $\mat{F}$. The $m\times 1$ vector $\vect{\Delta W}^{n+1/2}$ is a vector of independently distributed standard normal samples for white noise excitation. Specialized stochastic algorithm are necessary for white noise as energy is spread out in the frequency domain, so it is not possible to directly use a small time step to resolve responses.
The Euler-Maruyama scheme is used in \Cref{sec:2dof_example}.

For the aeroelastic and frictional systems, the stochastically varying parameters are perturbed, and the entire time step of $[t^n, t^{n+1}]$ is evaluated using the perturbed parameter vector with \eqref{eq:2ndOrderTimeRes} and \eqref{eq:1stOrderTimeRes}.
Consistent with the Euler-Maruyama scheme \cite{kloeden_numerical_1999}, a parameter is perturbed from the baseline value $p_i$ for the current time step as
\begin{equation}
    p_i^{n+1/2} = p_i + \dfrac{\sigma}{\sqrt{\Delta t}} \left[\vect{\Delta W}^{n+1/2}\right]_i,
\end{equation}
where the scalar $\left[\vect{\Delta W}^{n+1/2}\right]_i$ is sampled independently at each time step from the standard normal distribution.

\subsection{Harmonic Balance Method} \label{sec:appendix_methods_hbm}

HBM represents the periodic steady-state solution as a combination of harmonic components for each DOF \cite{krackHarmonicBalanceNonlinear2019}. For the frictional system, a standard HBM implementation from \cite{porter_efficient_2025, porterTMDSimPy} is used with minimal modifications noted in \Cref{sec:friction_system}. For this approach, the first-order frictional states are not explicitly variables in the HBM solution, but are calculated through the AFT procedure \cite{krackHarmonicBalanceNonlinear2019}.

For the aeroelastic system, HBM is implemented as a nonintrusive wrapper to the state-space model. In this case, HBM represents both the second-order states $\vect{y}$ and first-order states $\vect{z}$ as a series of harmonics.
Given harmonic coefficients for the current iteration, the time series representation of the states are calculated as $\vect{\tilde{y}}$ and $\vect{\tilde{z}}$ via the AFT procedure using an inverse of a fast Fourier transform to evaluate the responses at $2^n$ times in the interval $[0, T-dt]$ where $dt = T/2^n$. Then a time-varying residual is constructed at these time instants as
\begin{equation}
    \vect{R}(t) = 
    \begin{bmatrix}
        \vect{\ddot{\tilde{y}}}(t) \\ \vect{\dot{\tilde{z}}}(t)
    \end{bmatrix}
    -
    \begin{bmatrix}
        \mat{0} & \mat{I_y} & \mat{0} \\
        \mat{0} & \mat{0} & \mat{I_z}
    \end{bmatrix}
    \vect{f}\left( t, 
    \begin{bmatrix}
        \vect{\tilde{y}}(t) \\ \vect{\dot{\tilde{y}}}(t) \\ \vect{\tilde{z}}(t)
    \end{bmatrix},
    \vect{p} \right),
\end{equation}
where $\mat{I_y}$ and $\mat{I_z}$ are identity matrices of the sizes that correspond to $\vect{y}$ and $\vect{z}$, respectively. Here, $\mat{0}$ represents zero matrices of the appropriate sizes.
Like in the AFT procedure, the fast Fourier transform of the time-varying residual is evaluated to give a harmonic residual vector that is used in the HBM solution.

\section{Numerical Conditioning for Self-Excited Calculations}
\label{sec:app_conditioned_self_excited}

This work does not directly consider self-excited systems. However, this section presents an outline for modifying the self-excited calculation from \cite{ahsan_adjoint_based_2025} to ensure numerical conditioning consistent with \Cref{sec:updated_alg}.
For the self-excited case, $\mat{\mathcal{J}}$ is calculated in continuous form as
\begin{equation}
\begin{split}
    \mat{\mathcal{J}}(t^i) 
    &= \mat{Q}(t^i) \mat{X}(t^i) \int_0^{t^i} \mat{G}(s) \mat{G}^T(s) ds \mat{X}^T({t^i}) \mat{Q}^T(t^i),
\end{split}
\end{equation}
or in the discrete notation used here as
\begin{equation}
    \mat{\mathcal{J}}^i 
    = \mat{Q}^i \mat{X}^i \int_0^{t^i} \mat{G}(s) \mat{G}^T(s) ds \left[\mat{X}^i\right]^T \left[\mat{Q}^i\right]^T .
\end{equation}
Here,
\begin{equation}
    \mat{Q}^i = \mat{X}^i \left[ \mat{I_N} - \vect{\dot{x}}^0 \vect{w}^T \right] \left[\mat{X}^i \right]^{-1}, 
    \end{equation}
where $\vect{\dot{x}}^0$ is the time derivative of the states at time 0.
The vector $\vect{w}$ satisfies \cite{ahsan_adjoint_based_2025}
\begin{equation}
\begin{split}
    \vect{w}^T \vect{\dot{x}}^0 &= 1
    \\
    \vect{w}^T \mat{X}^{N_t} &= \vect{w}^T.
\end{split}
\end{equation}
Using the latter property, $\mat{Q}^i$ can be rewritten as
\begin{equation}
    \mat{Q}^i = 
    \mat{I_N}
    -
    \mat{X}^i \vect{\dot{x}}^0 \vect{w}^T \left[ \mat{X}^{N_t} \left(\mat{X}^i \right)^{-1}  \right],
    \end{equation}
where the term in the square brackets can now be calculated in a conditioned way rather than via inverting $\mat{X}^i$.
This numerically conditioned $\mat{Q}^i$ can then be used to pre- and post-multiply the $\mat{\mathcal{J}}^i$ that is calculated for the externally excited system in \Cref{sec:updated_alg}.

\contrib Justin H. Porter developed ideas and theory, implemented numerical solutions, and wrote manuscript.

\ack

This work was authored by the National Laboratory of the Rockies for the U.S. Department of Energy (DOE), operated under Contract No. DE-AC36-08GO28308. This work was supported by the Laboratory Directed Research and Development (LDRD) Program at the National Laboratory of the Rockies. The views expressed in the article do not necessarily represent the views of the DOE or the U.S. Government. The U.S. Government retains and the publisher, by accepting the article for publication, acknowledges that the U.S. Government retains a nonexclusive, paid-up, irrevocable, worldwide license to publish or reproduce the published form of this work, or allow others to do so, for U.S. Government purposes.

\printbibliography{yourbib}

\end{document}